\documentclass[a4paper,10pt]{article}

\usepackage[top=4cm, bottom=3cm, left=3.2cm, right=3.2cm]{geometry}

\usepackage{Macros}

\begin{document}

\title{Continuous approximation of $M_t/M_t/ 1$ distributions with application to production}

\author{Dieter Armbruster\footnotemark[1], \; Simone G\"ottlich\footnotemark[2], \; Stephan Knapp\footnotemark[2]}

\footnotetext[2]{Arizona State University, School of Mathematical and Statistical Sciences, Tempe, AZ 85257-1804 (armbruster@asu.edu)}
\footnotetext[2]{University of Mannheim, Department of Mathematics, 68131 Mannheim, Germany (goettlich@uni-mannheim.de, stknapp@mail.uni-mannheim.de).}

\date{\today}

\maketitle

\begin{abstract}
A single queueing system with time-dependent exponentially distributed arrival processes and exponential machine processes (Kendall notation $M_t/M_t/1$) 
is analyzed. Modeling the time evolution for the discrete queue-length distribution by a  continuous drift-diffusion process a Smoluchowski equation on the half space  is derived 
approximating the forward Kolmogorov equations. The approximate model is analyzed and validated, showing excellent agreement for the probabilities of all queue lengths and for all  queuing utilizations,  including ones that are very small and some that are significantly larger than one. Having an excellent approximation for the probability of an empty queue generates an 
approximation of the expected outflow of the queueing system. Comparisons to several well-established approximation from the literature show significant improvements in several numerical examples.

\end{abstract}


\maketitle

\section{Introduction}
We consider the case of a single queue with one server,  a FIFO (first-in first-out) service rule and  Markovian arrival and departure processes with a time-dependent distribution. This reads in Kendall's notation $M_t/M_t/1$.
Thus the system is specified by random arrival of goods and random service times, where the distribution is dependent on time but the system as a whole is Markovian. Determining the expected outflow of the production is an important performance measure and is based on the idle probability of the queuing system. In the stationary case, the results are well-known but the transient case is difficult to tackle and a system of infinitely many ordinary differential equations (ODEs) has to be solved. 
Specifically the transient behavior of the probability $p_k(t)$ having $k$ customers in the queueing system (queue and service) at time $t \geq 0$ is for time-independent rates given by a formula consisting of an infinite sum of modified Bessel functions of the first kind; see \cite{Gross_Queueing}.

Our interest in this problem comes from the concept of production planning  and control in manufacturing industries which has been around since the early 20th century
and has a long history of research (for reviews see \cite{armbruster2012continuous, kempf2011planning}).  In essence the problem is to determine the input into a production resource  to generate
a desired output over time. Operationally the problem splits into a  forward problem, which  estimates the expected output trajectory given a specified input trajectory, and a  backward problem which determines the input pattern required to produce a desired output pattern on average over time. 

As production starts,  availability of parts, machines and workers are all determined by random processes, there is a long tradition to discuss production planning in the context of queueing theory. In particular performance measures such as the average time in system or queue length under long-run steady-state conditions can be derived for e.g.\ queuing networks of Jackson type \cite{chen2013fundamentals}. However for many manufacturing processes, notably in the semi-conductor industry, cycle times are long and planning period short relative to the cycle time violating the assumptions of a steady-state approximation  and thus leading to the analysis of time-dependent queuing networks. 

At the same time the backward problem is an optimization or an optimal control problem: given a desired output trajectory, find the optimal input function under the constraint that input and output are related via the solution of the forward problem.  If the forward problem can be described by an evolution equation (a set of ordinary
or partial differential equations), then such problems can be solved using adjoint calculus \cite{la2010control}. In particular, existence and uniqueness  as well as controllability  of solutions  can be proven \cite{ colombo2011control, coron2012controllability,keimer2017existence}
in some cases. 

This suggests an attempt to model a time-dependent queueing system via a continuous description.  There have been two different strands of research in this direction 
in the last 50 years: Newell in a series of papers in the 1960 suggested a {\em diffusion approximation} \cite{Newell11986,Newell21968, Newell31968} and postulated a Fokker Planck equation for the 
cumulative  distribution function $F(x,t)$  for the queue length $x$ at time $t$. He created models for traffic flows though rush hour. 
A different approach was introduced by \cite{ArmbrusterMarthalerRinghofer, DegondRinghofer2007} based on kinetic theory for the probability density $f(x,v,t)$ of finding a particle at position $x$ in the production process considered as a  queue and a machine, moving forward with speed $v$ at time $t$.  Boltzmann equations  for $f$ and  moment equations with different closures lead to  transport equations for the density  similar to the Lighthill-Whitham-Roberts model \cite{lighthill1955kinematic,richards1956shock}
  for traffic flow and to second order moment equations for the velocity of particles moving through the queues \cite{ArmbrusterMarthalerRinghofer}.

Recently, Armbruster et al.\ \cite{armbrusterkinetic} performed a systematic analysis comparing simulations of these moment equations with discrete event simulations (DES) 
for a  factory production  modeled as an $M/M/1$ queue with a non-homogeneous Poisson arrival process. They show  that while using higher order moment equations 
improves the model, these transport equations have  intrinsic timescales  that are not present in the original stochastic processes and thus may lead to fundamentally bad approximations in some cases. 

Our research picks up Newell's idea of a Fokker-Planck equation for the queue length. We derive an approximation for the probability $p_k(t)$ of
having  $k \in \N_0$ customers in the system at time  $t$  through a continuous variable $\rho(x,t)$ leading to a drift-diffusion equation known as the Smoluchowski equation. 
We study the relationship of the queue-length probabilities and the solution to the Smoluchowski equation on the half space $[0,\infty)$ with a linear potential. Additionally, we can find the explicit solution to this equation in the transient case with time-independent rates, based on calculations by  Smoluchowksi \cite{Smoluchowski1916Vortrag,Smoluchowski1916Paper}.  
Additionally, in the case of time-dependent rates, a numerical scheme is provided 
and compared to the ODE system solution as well as approximate formulas taken from \cite{MasseyPender2013, Rider1976}. 
%
%
\subsubsection*{Classification of the Proposed Continuous Approach}

There have been a number of previous studies generating various approximations for time-varying queues. Whitt \cite{Whitt2018} in a recent review presents the historical development of methods and discusses the relevant literature.
The first and simplest approach is the so-called {\em Pointwise Stationary Fluid Flow Approximation} introduced by Rider \cite{Rider1976} and later used  in \cite{Wang1996} leading to ODEs for the expected queue length. We discuss some details of this in section 6. 

Mandelbaum \cite{Mandelbaum1995} uses a first order macroscopic scaling of influx rate and production rate $\lambda \Rightarrow \lambda^{\epsilon} = \frac{\lambda}{\epsilon}$, 
$\mu \Rightarrow \mu^{\epsilon} = \frac{\mu}{\epsilon}$ proving the asymptotic validity of a fluid approximation for the queue length $L^\epsilon$. 
The latter approximation has been improved by a mesoscopic diffusion approximation and analyzed for different traffic regimes.
These approximations can be used for a local asymptotic expansion of the queue length and again assume stationarity in the limit. Since this is similar to \cite{Rider1976, Wang1996}, we omit a comparison later on.

In \cite{MasseyPender2013}, the \textit{fluid mean} (first order), \textit{Gaussian variance} (second order) and \textit{Gaussian skewness approximation} (third order) are introduced and compared. In particular, the Gaussian variance and skewness approximation lead to fairly good results and we compare our approach to these methods in section 6. 

There is also a so-called Poisson-Charlier approximation, see \cite{Pender2014} in which the queue length process is approximated by a truncated series of Poisson-Charlies polynomials.
The first order expansion leads exactly to the fluid mean approximation from above and only performs well in the case of a large numbers of servers, which is not the case in the model we consider (one server). A second order expansion is also introduced but we expect a better performance of the Gaussian skewness approximation for the $M_t/M_t/1$ queueing model.

This paper is organized as follows. In the second section we state basic definitions and results from the queueing theory, which is needed throughout this paper. The third section addresses the formal derivation of the approximate model for an $M_t/M_t/1$ queueing distribution, followed by explicit solutions to the associated Smoluchowski equation. This section is followed by the numerical treatment of this model, which is used to study the approximate model and compare it to the exact solution in various examples in section 5. In section 6, the connection to the production context is introduced and the approximate model is compared to several well-established approximation from the literature in numerical examples.

\section{Definition and Basic Results on $M_t/M_t/1$ Queues}
 Let $\lambda \colon \R_{\geq 0} \to \R_{\geq 0}$ denote the time-dependent arrival rate, and let $\mu  \colon \R_{\geq 0} \to \R_{\geq 0}$ denote the time-dependent processing rate. We denote the number of customers in the system at time $t \geq 0$ with the random variable $L(t)$ and define $p_k(t) = P(L(t) = k)$ as the probability to have $k \in \N_0$ customers in the system at time $t$.
Specifically, we study the behavior of the queueing system given by figure \ref{fig:GraphSingleQueueServerTime1}.
 \definecolor{mygrey}{RGB}{100,100,100}
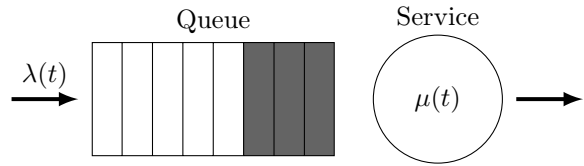
\begin{figure}[H]
\centering
\begin{tikzpicture}[>=latex]
\node[
  myshape,
  rectangle split part fill={white,white,white,white,white,mygrey}
  ] 
  (shape) {};
\node[draw,circle,inner sep=0.6cm,right=of shape,shift={(-5mm,-0mm)}]
  (out) {};

\draw[line width=1.5pt,->] 
    ([xshift=5pt]out.east|-out.east) -- ([xshift=30pt]out.east|-out.east);
\draw[line width=1.5pt,->] 
    ([xshift=-30pt]shape.west|-shape.west) to node[above]{$\lambda(t)$} ([xshift=-5pt]shape.west|-shape.west);

\node[anchor=north,above] 
  at (shape|-shape.north)
  {Queue};
\node[anchor=north,above] 
  at (out|-out.north)
  {Service};
\node
	at (out|-out)
	{$\mu(t)$};
\end{tikzpicture}
\caption{Graphical representation of a single queue service unit}
\label{fig:GraphSingleQueueServerTime1}
\end{figure}
 
Since the queue length follows a birth-death process, we formally obtain the following relation for $k\geq 1$:
\begin{align*}
p_k(t+\Delta t) &= P(L(t+\Delta t) = k)\\
&= P(L(t+\Delta t) = k|L(t) = k)p_k(t) \\
&+ P(L(t+\Delta t) = k|L(t) = k-1)p_{k-1}(t) \\
&+ P(L(t+\Delta t) = k|L(t) = k+1)p_{k+1}(t) +o(\Delta t)\\
&= (1-\Delta t(\lambda(t)+\mu(t)))p_k(t) + \Delta t \lambda(t) p_{k-1}(t) +\Delta t \mu(t) p_{k+1}(t) + o(\Delta t),
\end{align*}
leading in the limit $\Delta t \to 0$ to the Kolmogorov forward equations
\begin{align}
\frac{d}{dt}p_k(t) = \lambda(t) p_{k-1}(t)-(\lambda(t)+\mu(t))p_{k}(t)+\mu(t) p_{k+1}(t). \label{eq:KFTransientMM11}
\end{align}
For $k = 0$, we obtain
\begin{align}
\frac{d}{dt}p_0(t) = \mu(t)p_1(t)-\lambda(t)p_0(t). \label{eq:KFTransientMM12}
\end{align}
 
The situation of time-independent arrival and processing rates is well established and can be found in the standard literature, e.g., \cite{Gross_Queueing}. The steady-state distribution in the case 
$\varrho := \frac{\lambda}{\mu}<1$  is given by 
\begin{align}
\overline{p}_k = (1-\varrho)\varrho^k, \label{eq:SS_pdf}
\end{align} 
where $\varrho$ is the so-called traffic intensity or utilization of the queueing system. Since the steady-state distribution is a geometric distribution, the expected queue length is given by \[\E [\overline{L}] = \lim_{t \to \infty} \E [L(t)] = \frac{\varrho}{1-\varrho}.\] 
 
\section{Continuous Approximation of $M_t/M_t/1$ Distributions}
%
\subsection{Derivation of the Approximating Model}
Instead of approximating the exact solution $p_k(t)$ to eq. (\ref{eq:KFTransientMM11}), we derive an approximate model that we solve exactly (or numerically). In \cite{Newell11986,Newell21968,Newell31968}, the cumulative distribution function $F(t,x)$ of $L(t)$ is approximated by a second-order PDE. We follow this idea, but we model a probability density function (pdf) instead. Specifically, let us assume that $\rho \colon \R_{\geq 0}^2 \to \R_{\geq 0}$ is the solution to
\begin{alignat}{2}
\rho_t(x,t) +a(t) \rho_x(x,t)-b(t) \rho_{xx}(x,t) &= 0 &&\quad\text{ on } \R_{>0}^2,\label{eq:AppQueueModelEq1}\\
a(t) \rho(x,t) - b(t) \rho_x(x,t) &= 0 &&\quad\text{ on } \{0\}\times\R_{\geq 0},\label{eq:AppQueueModelEq2}\\
\rho(x,t) &= \rho_0(x) &&\quad\text{ on } \R_{\geq 0} \times \{0\},\label{eq:AppQueueModelEq3}
\end{alignat}
and we want to derive conditions on the functions $a$, $b$ and $\rho_0$ such that
\begin{align*}
\int_k^{k+1} \rho(x,t) dx \approx p_k(t).
\end{align*}
The coefficient $a(t)$ will describe the mean behavior of the model, i.e., $\gl(t)>\mu(t)$ implies an increasing size of customers such that the probability density function is expected to move to the right; in the case  $\gl(t)<\mu(t)$, it is exactly the opposite. Since the system is not deterministic, the coefficient $b(t)$ inherits the variance or fluctuations of the model.
Assuming  $\lim_{x \to \infty }\rho(x,t) = \lim_{x \to \infty }\rho_x(x,t) = 0$ we observe  conservation of  mass:
\begin{align*}
\frac{d}{dt} \int_0^\infty \rho(x,t) dx =  \int_0^\infty -a(t) \rho_x(x,t)+b(t) \rho_{xx}(x,t) dx = a(t) \rho(0,t)-b(t)\rho_x(0,t) = 0,
\end{align*}
Since we consider a pdf we set
\begin{assumption}\label{ass:unit_mass}
\[\int_0^\infty \rho_0(x) dx = 1.\]
\end{assumption}

We compute the change in the expected number of customers using  \eqref{eq:KFTransientMM11} as
\begin{align*}
\frac{d}{dt} \E [L(t)] &= \sum_{k = 0}^\infty k \frac{d}{dt} p_k(t)\\
& =\sum_{k=1}^\infty k  \lambda(t) p_{k-1}(t)-(\lambda(t)+\mu(t))p_{k}(t)+\mu(t) p_{k+1}(t).
\end{align*}
Rearranging the terms yields
\begin{align}
\frac{d}{dt} \E [L(t)]   =&\;\sum_{k=1}^\infty k \left( \lambda(t) (p_{k-1}(t)-p_k(t))+\mu(t) (p_{k+1}(t)-p_k(t)) \right)\notag\\
 =&\; \gl(t) \left(\sum_{k=1}^\infty k p_{k-1}(t) - \sum_{k = 1}^\infty k p_k(t) \right)\notag\\
 &\;+\mu(t) \left(\sum_{k = 1}^\infty k p_{k+1}(t)- \sum_{k = 1}^\infty k p_k(t) \right)\notag\\
  =&\; \gl(t) \left(\sum_{k=0}^\infty (k+1) p_{k}(t) - \E[L(t)] \right)\notag\\
 &\;+\mu(t) \left(\sum_{k = 2}^\infty (k-1) p_{k}(t)- \E[L(t)] \right)\notag\\[1ex]
   =&\; \gl(t) \left(\E[L(t)]+ 1 - \E[L(t)] \right)\notag\\[1ex]
 &\;+\mu(t) \left(\E[L(t)]-p_1(t)-(1-(p_0(t)+p_1(t)))- \E[L(t)] \right)\notag\\[1ex]
=&\; \lambda(t)-\mu(t) + \mu(t) p_0(t) \label{eq:queueingExpectation};
\end{align}
see, e.g., \cite{Rider1976}. This implies that the rate of change of the expected number of customers in the system is given by the arrival minus the service rate plus the service rate multiplied with the idle probability $p_0(t)$. \\
On the other hand, integration by parts, assuming $\lim_{x \to \infty }x\rho(x,t) = \lim_{x \to \infty }x\rho_x(x,t) = 0,$ implies
\begin{align}
\frac{d}{dt} \int_0^\infty x \rho(x,t)dx &=-a(t)\int_0^\infty x \rho_x(x,t) dx + b(t) \int_0^\infty x \rho_{xx}(x,t) dx \notag \\
& = a(t) \int_0^\infty \rho(x,t) dx - b(t) \int_0^\infty \rho_x(x,t) dx \notag \\[1ex]
&= a(t)+b(t) \rho(0,t). \label{eq:PDEExpecation}
\end{align}
Comparing \eqref{eq:queueingExpectation} and \eqref{eq:PDEExpecation}, leads to  
\begin{assumption}\label{ass:fun_a}
\[a(t) = \gl(t)-\mu(t)\]
\end{assumption}
\noindent to be satisfied at every time $t \geq 0$.\\

If $\varrho(t) := \frac{\lambda(t)}{\mu(t)} \geq 1$, we have $p_0(0,t) \approx \rho(0,t) \approx 0$, and assumption \ref{ass:fun_a} corresponds to the expected increase in the number of customers in the system.
Assumption \ref{ass:fun_a} guarantees the inclusion of the mean transient behavior.  Additionally we want the approximate model to be exact in the steady state ($\varrho<1$) such that eventually occurring errors become damped. This implies that
\begin{align*}
\overline{p}_k = \int_k^{k+1} \overline{\rho}(x) dx,
\end{align*}
where $ \overline{\rho}$ is the steady-state solution to \eqref{eq:AppQueueModelEq1}-\eqref{eq:AppQueueModelEq3} and $\overline{p}_k$ is the steady state
distribution given in \eqref{eq:SS_pdf}. The steady state solution to \eqref{eq:AppQueueModelEq1} is given by 
$ \overline{\rho}(x) =  \tilde{C} e^{\frac{a}{b}x}.$  Since $\rho(\cdot,t)$ and $\overline{\rho}$ are probability density functions, we have 
 \begin{align}
 \overline{\rho}(x)  = -\frac{a}{b}e^{\frac{a}{b}x} \label{eq:SteadyStateAppQueue}.
 \end{align}
 Thus
$\int_k^{k+1} \overline{\rho}(x)dx = e^{\frac{a}{b}k}(1-e^{\frac{a}{b}}),$
 which equals $\overline{p}_k$ if and only if $\frac{a}{b} = \ln(\varrho) = \ln(\lambda)-\ln(\mu)$, leading to 
 \begin{assumption}
  \[b(t) = \frac{\mu(t)-\lambda(t)}{\ln(\mu(t))-\ln(\lambda(t))}.\]
\end{assumption} 
 Altogether, 
\begin{align*}
p^A_k(t) = \int_k^{k+1} \rho(x,t) dx,
\end{align*}
is an approximation of the distribution of the queueing model $p_k(t)$ 
where $\rho$ is the solution to \eqref{eq:AppQueueModelEq1}-\eqref{eq:AppQueueModelEq3} with the coefficients
\begin{align}
a(t) &= \gl(t)-\mu(t), \label{eq:QueueingKoeffA}\\
b(t) &= \frac{\mu(t)-\lambda(t)}{\ln(\mu(t))-\ln(\lambda(t))},\label{eq:QueueingKoeffB}
\end{align}
and initial data $\rho_0$ satisfying
$\int_k^{k+1} \rho_0(x)dx = p_k(0)$ for every $k \in \N_0$.

\subsection{Analytic Solution}
In the time-homogeneous case, $\lambda(t) \equiv \lambda$ and $\mu(t) \equiv \mu$, equations \eqref{eq:AppQueueModelEq1}-\eqref{eq:AppQueueModelEq2} are a simple case of the so-called Smoluchowski equations on the half space $[0, \infty)$; see \cite{Lamm1983, Smoluchowski1916Vortrag, Smoluchowski1916Paper}. Smoluchowski derived the fundamental solution $\rho^F$ to this equation, which in our context reads
\begin{align}
\rho^F(x,t|x_0,0) = \rho^{F,1}(x,t|x_0,0)+\rho^{F,2}(x,t|x_0,0)+\rho^{F,3}(x,t|x_0,0)
\end{align}
for $x_0,x,t \in \R_{\geq 0}$
with
\begin{align*}
\rho^{F,1}(x,t|x_0,0) &= \frac{1}{\sqrt{4 \pi b t}}e^{-\frac{(x-x_0-a t)^2}{4 b t}},\\
\rho^{F,2}(x,t|x_0,0) &= \frac{1}{\sqrt{4 \pi b t}}e^{-\frac{a}{b}x_0-\frac{(x+x_0-a t)^2}{4 b t}},\\
\rho^{F,3}(x,t|x_0,0) &= -\frac{a}{2 b}e^{\frac{a}{b}x} \frac{2}{\sqrt{\pi}} \int_{\frac{x+x_0+at}{\sqrt{4 b t}}}^\infty e^{-y^2}dy.
\end{align*}
Analogously to the heat equation, the function
\begin{align}
\rho(x,t) = \int_0^\infty \rho^F(x,t|x_0,0) \rho_0(x_0) dx_0
\end{align}
solves the initial boundary value problem, see \eqref{eq:AppQueueModelEq1}-\eqref{eq:AppQueueModelEq3}, provided $\lambda(t) \equiv \lambda$, $\mu(t) \equiv \mu$, see \cite{StraussPDE2008}.

Since we are interested in the transient behavior of the queueing distribution, we assume that we  start with a steady-state distribution, see \eqref{eq:SteadyStateAppQueue}, which is in the form of
\begin{align}
\rho_0(x) = \overline{\rho}(x,0) = -\frac{a_0}{b_0}e^{\frac{a_0}{b_0}x}, \label{eq:rho0StSt}
\end{align}
where $a_0 = \lambda_0-\mu_0$, $b_0 = \frac{\mu_0-\lambda_0}{\ln(\mu_0)-\ln(\lambda_0)}$, as in \eqref{eq:QueueingKoeffA}-\eqref{eq:QueueingKoeffB}, and provided $\frac{\lambda_0}{\mu_0}<1$. We compare the transition to the new steady state given by
\begin{align*}
\lim_{t \to \infty}\rho(x,t) = -\frac{a}{b}e^{\frac{a}{b}x}
\end{align*}
determined by $\lambda$ and $\mu$, again provided $\frac{\lambda}{\mu}<1$.
In this case, we can derive the solution to \eqref{eq:AppQueueModelEq1}-\eqref{eq:AppQueueModelEq3} with \eqref{eq:rho0StSt} explicitly; it reads
\begin{align}
\rho(x,t) =&\; \rho^1(x,t)+\rho^2(x,t)+\rho^3(x,t),\label{eq:rhoStSt1}\\[1ex]
\rho^1(x,t) =&\; -c_0 e^{c_0(c_0 b t+x-a t)} \Phi\left(\frac{c_02bt+x-at}{\sigma} \right),\notag \\
\rho^2(x,t) =&\; -c_0 e^{d(d b t-x+at)} \Phi\left(\frac{d 2bt-x+at}{\sigma}\right),\notag \\
\rho^3(x,t) =&\; -c e^{c x} \Phi \left(- \frac{x+at}{\sigma}\right)\notag 
+c e^{c x} e^{c_0(c_0 b t-x-a t)} \Phi \left(\frac{c_0 2 b t-x-at}{\sigma}\right),
\end{align}
with 
\begin{alignat*}{4}
c &= \frac{a}{b},\quad c_0 &&= \frac{a_0}{b_0},\quad d &&= c_0-c, \quad \sigma &&= \sqrt{2bt} \quad \text{and} \quad \Phi(z) = \int_{-\infty}^z \frac{1}{\sqrt{2 \pi}}e^{-\frac{x^2}{2}}dx.
\end{alignat*}
To calculate the approximated probabilities $p_k^A(t)$ of $p_k(t)$, we have to integrate \eqref{eq:rhoStSt1} over $[k,k+1)$. A calculation results in
\begin{align}
p_k^A(t) =&\; \frac{c_0-c}{d}\left(\Phi\left(\frac{k+1-at}{\sigma}\right)-\Phi\left(\frac{k-at}{\sigma}\right)\right)\notag \\
&\;-e^{c_0(c_0b-a)t}e^{c_0k}\left(e^{c_0}\Phi\left(\frac{(2c_0b-a)t + k+1}{\sigma}\right)-\Phi\left(\frac{(2c_0b-a)t+k}{\sigma}\right)\right)\notag\\
&\;+\frac{c_0}{d} e^{d(db+a)t}e^{-dk}\left(e^{-d}\Phi\left(\frac{(2db+a)t - (k+1)}{\sigma}\right)-\Phi\left(\frac{(2db+a)t-k}{\sigma}\right)\right)\notag\\
&\;-e^{ck}\left(e^{c}\Phi\left(-\frac{k+1+at}{\sigma}\right)-\Phi\left(-\frac{k+at}{\sigma}\right)\right)\notag\\
&\;-\frac{c}{d} e^{c_0(c_0 b-a)t}e^{-dk}\left(e^{-d}\Phi\left(\frac{(2c_0b-a)t - (k+1)}{\sigma}\right)-\Phi\left(\frac{(2c_0b-a)t-k}{\sigma}\right)\right).  \label{eq:PStSt1}
\end{align}

We additionally discuss an alternative way to calculate the integrals in the following. Using the mean value theorem 
\begin{align}
\int_k^{k+1} \rho(x,t)dx = \rho(\xi_k(t),t),
\end{align}
we reduce the effort to one evaluation at $\xi_k(t)$ of $\rho$. If we can provide a good approximation of $\xi_k(t)$, we reduce the numerical costs in the time-dependent case as well. The continuous-model approximation is exact in steady state, with 
\begin{align*}
\overline{p}_k(t) = \int_{k}^{k+1} \overline{\rho}(x) dx = \varrho^k(1-\varrho) = -\ln(\varrho)e^{\ln(\varrho)\overline{\xi}_k},
\end{align*}
which is equivalent to \[\overline{\xi}_k = k+ \frac{\ln(1-\varrho)-\ln(-\ln(\varrho))}{\ln(\varrho)}.\] This motivates the use of 
\begin{align*}
\tilde{\xi}_k(t) = 
\begin{cases}
k & \text{ if } \varrho(t) = 0,\\
k+ \frac{\ln\left(\frac{\varrho(t)-1}{\ln(\varrho(t))}\right)}{\ln(\varrho(t))} &\text{ if } \varrho(t)\in (0,\infty)\setminus \{1\},\\[1ex]
k+\nicefrac{1}{2} &\text{ if } \varrho(t) = 1,
\end{cases}
\end{align*}
 which is continuous in $\varrho(t)$ and satisfies $\tilde{\xi}_k(t) \in [k,k+1]$. We define the approximation
\begin{align}
\tilde{p}_k(t) = \rho(\xi_k(t),t). \label{eq:PXi1}
\end{align}

\section{Numerical Treatment}
To compare and validate the continuous approximation \eqref{eq:AppQueueModelEq1}-\eqref{eq:AppQueueModelEq3} with the result of the system of ordinary differential equations \eqref{eq:KFTransientMM11}-\eqref{eq:KFTransientMM12}, we need an approximation of $p^A_k(t)$ and $p_k(t)$. The latter is approximated by reducing the infinite ODE system to $N \in \N$ equations, where a ``boundary'' condition is set, i.e., we use the mass conservation to close the equations. We have
\begin{align*}
\frac{d}{dt} \sum_{k=0}^{N-1} p_k(t) = \mu(t) p_N(t)-\lambda(t) p_{N-1}(t) \stackrel{!}{=} 0
\end{align*}
and set $p_N(t) = \frac{\lambda(t)}{\mu(t)}p_{N-1}(t)$
such that
\begin{align}
\frac{d}{dt} p_{N-1}(t) = \gl(t)p_{N-2}(t)-\mu(t) p_{N-1}(t). \label{eq:RestrictODESystem}
\end{align}
The resulting system is solved with the Matlab ODE solver \texttt{ode23}\footnote{Documentation: \href{https://de.mathworks.com/help/matlab/ref/ode23.html}{https://de.mathworks.com/help/matlab/ref/ode23.html}}.

To approximate $p_k^A(t)$, we distinguish two cases. In the case of constant coefficients $a,b$ and starting with a steady-state distribution given by $a_0,b_0$, we can use the analytic formulas \eqref{eq:rhoStSt1} and \eqref{eq:PStSt1} in combination with \eqref{eq:PXi1}. The evaluation of the standard normal cumulative distribution function is performed with the Matlab function \texttt{normcdf}\footnote{Documentation: \href{https://de.mathworks.com/help/stats/normcdf.html}{https://de.mathworks.com/help/stats/normcdf.html}}.

In the case of time-dependent coefficients $a,b$, we impose in the following a numerical scheme to approximate the solution of
\eqref{eq:AppQueueModelEq1}-\eqref{eq:AppQueueModelEq3}.
Let $\{i \Delta x \colon i \in \N_0\}$ be a discretization of the half space $[0,\infty)$ with fineness $\Delta x>0$. Since equations \eqref{eq:AppQueueModelEq1}-\eqref{eq:AppQueueModelEq3} imply the conservation of mass and since equation \eqref{eq:AppQueueModelEq2} is a no-flux boundary condition, it is natural to start with a conservative numerical scheme; see \cite{LeVequeRed}. Let $\rho_i^j$ be the approximation of $\rho(i \Delta x, t_j)$ for some time $t_j \geq 0$, and let
\begin{align*}
F(t,u,v) = -\frac{\mu(t)-\gl(t)}{\ln(\mu(t))-\ln(\gl(t))} \frac{(v-u)}{\Delta x}+(\gl(t)-\mu(t))\frac{u+v}{2}
\end{align*}
be the numerical flux function. We define the iteration by
\begin{align}
\rho_0^{j+1} &= \rho_0^j-\frac{\Delta t_j}{\Delta x}(F(t_j,\rho_0^j,\rho_{1}^j))-0),\label{eq:NumSchemeAppQueue1}\\
\rho_i^{j+1} &= \rho_i^j-\frac{\Delta t_j}{\Delta x}(F(t_j,\rho_i^j,\rho_{i+1}^j)-F(t_j,\rho_{i-1}^j,\rho_i^j)) \text{ for } i \in \N. \label{eq:NumSchemeAppQueue2}
\end{align}
In the case $\lambda(t_j)=\mu(t_j)$, which implies $a(t_j) = 0$, we only observe diffusion with $b(t_j) = \gl(t_j)$, and the solution should decrease in this time step. This is the case if we assume the standard stability condition for diffusion equations, see \cite{GRS}, which reads
\begin{align*}
\Delta t_j \leq \frac{(\Delta x)^2}{2}\frac{\ln(\mu(t_j))-\ln(\gl(t_j))}{\mu(t_j)-\gl(t_j)}
\end{align*}
to be satisfied in every iteration. We use the forward difference in time, which implies the resulting scheme to be first-order accurate in time, and from the second discrete and central derivative, we have second-order accuracy in space.

\section{Computational Results}
In the following, we numerically examine the continuous approximations $p^A_k(t)$ and $\tilde{p}_k(t)$ of the queue-length distribution $p_k(t)$. In the first part, we consider a steady state at time $t = 0$, which is determined by the rates $\gl_0>0$ and $\mu_0>0$ satisfying $\frac{\gl_0}{\mu_0}<1$, and the system has an abrupt change to the rates $\gl_1>0$ and $\mu_1>0$. In this case, we derived the analytic expression for $p^A_k(t)$ and $\tilde{p}_k(t)$, see \eqref{eq:PStSt1} and \eqref{eq:rhoStSt1}, such that a numerical scheme for the PDE is not necessary, and we avoid errors arising from the scheme. 

The second part addresses the use of the numerical scheme \eqref{eq:NumSchemeAppQueue1}-\eqref{eq:NumSchemeAppQueue2}, and we analyze the continuous approximation in the transient case, i.e., time-dependent rates. 

We introduce the following measures to evaluate the accuracy of the continuous approximations. Let $K \in \N$ be the number of equations that we want to compare, we define  the supremum error
\[ \|\epsilon\|_\infty = \max_{k \in \{ 0,\dots,K-1\}} |\epsilon_k|\] with $\epsilon \in \{\epsilon_k^A(t),\tilde{\epsilon}_k(t)\}$ defined by
\[\epsilon_k^A(t) = p_k(t)-p_k^A(t) \text{ and } \tilde{\epsilon}_k(t) = p_k(t)-\tilde{p}_k(t).\] 

We consider a time horizon $T = 100$ and $K = 100$ equations in all the cases, and we restrict the ODE system to $N = 1000$ equations; see \eqref{eq:RestrictODESystem}.

\subsubsection*{Ramp Up}
\begin{figure}[htb!]
\subfigure[Approximation of $p_0(t)$]{
\includegraphics[width=0.3\textwidth]{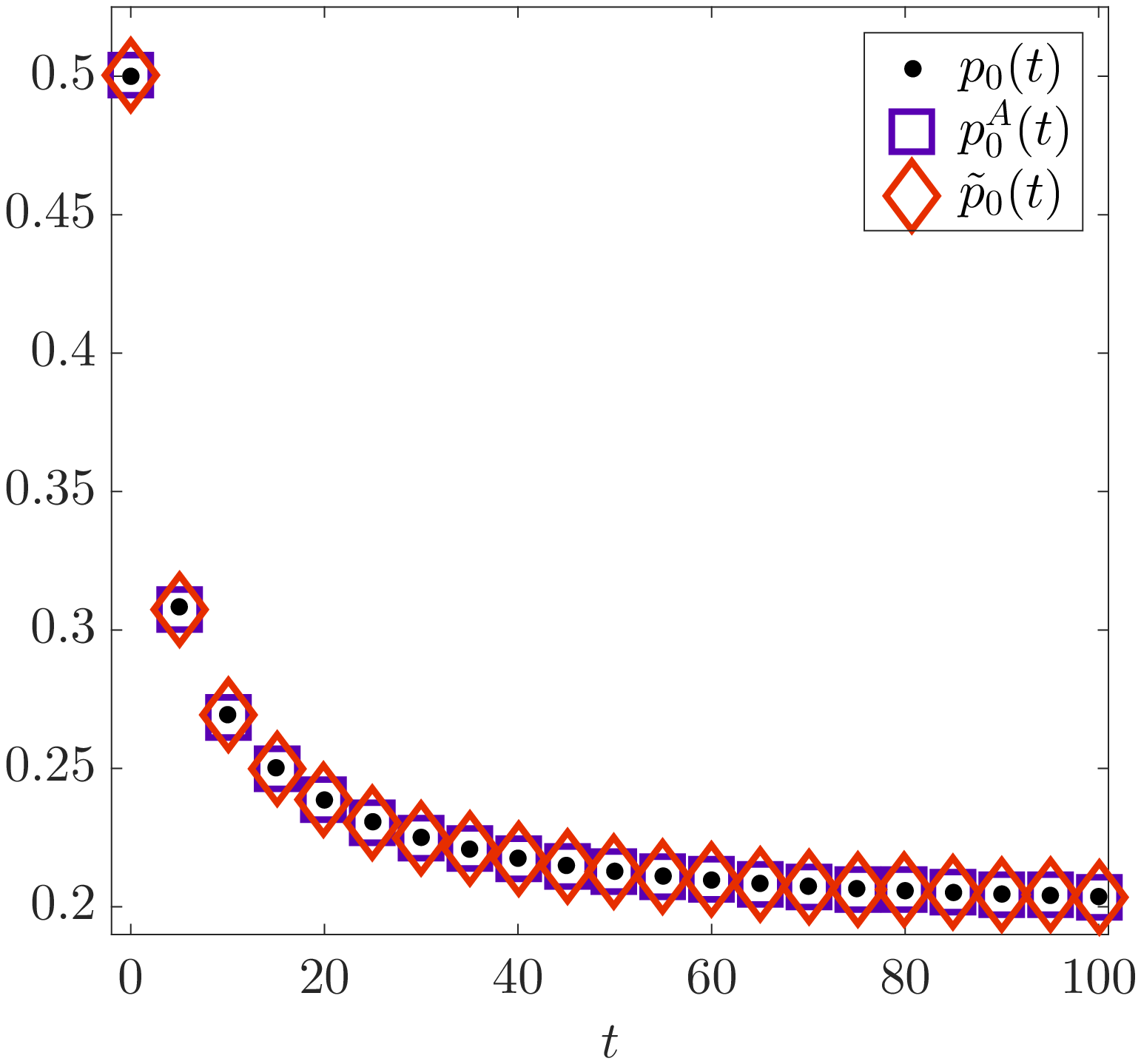}
}
\subfigure[Approximation of $p_3(t)$]{
\includegraphics[width=0.3\textwidth]{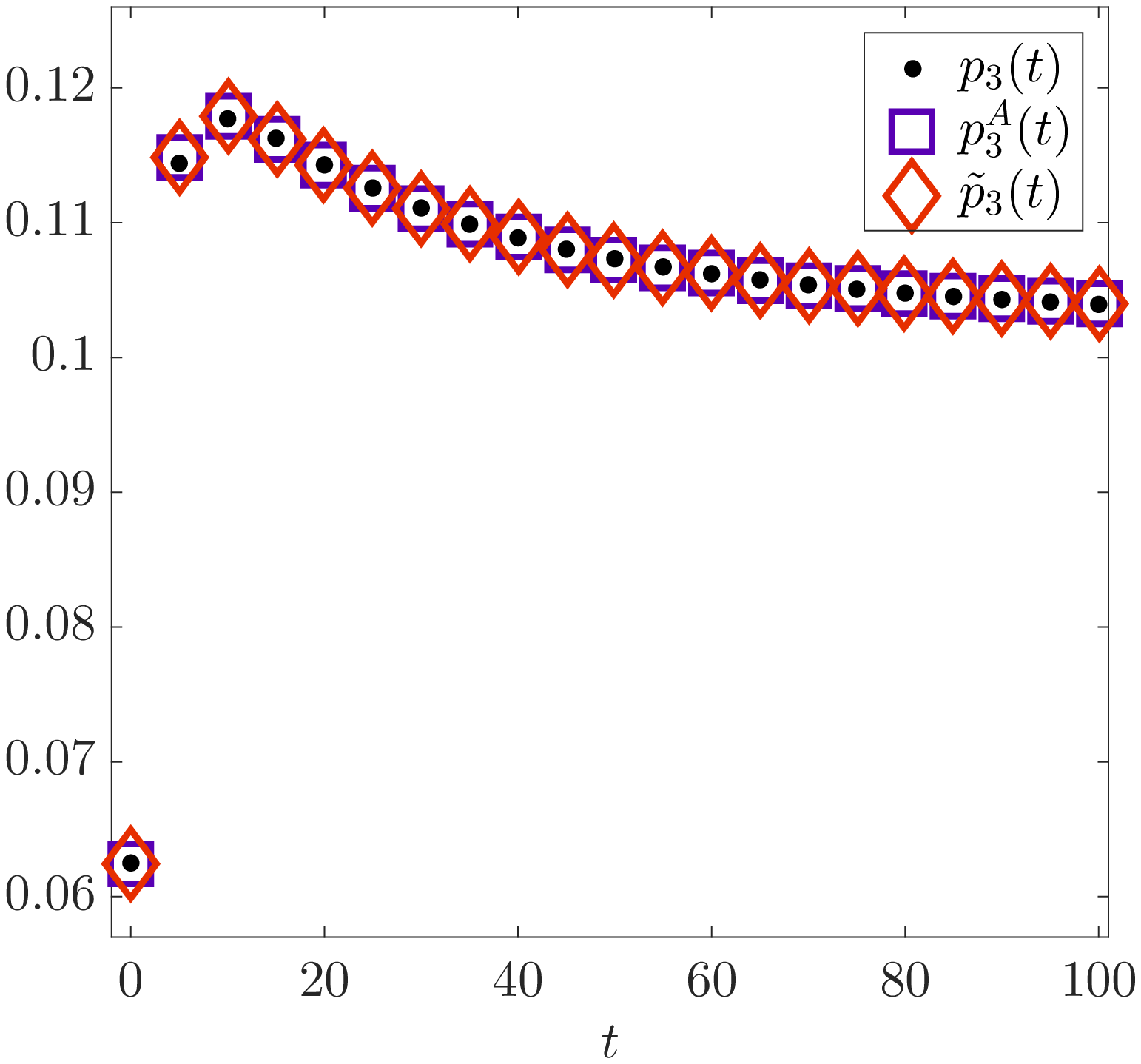}
}
\subfigure[Supremum error]{
\includegraphics[width=0.3\textwidth]{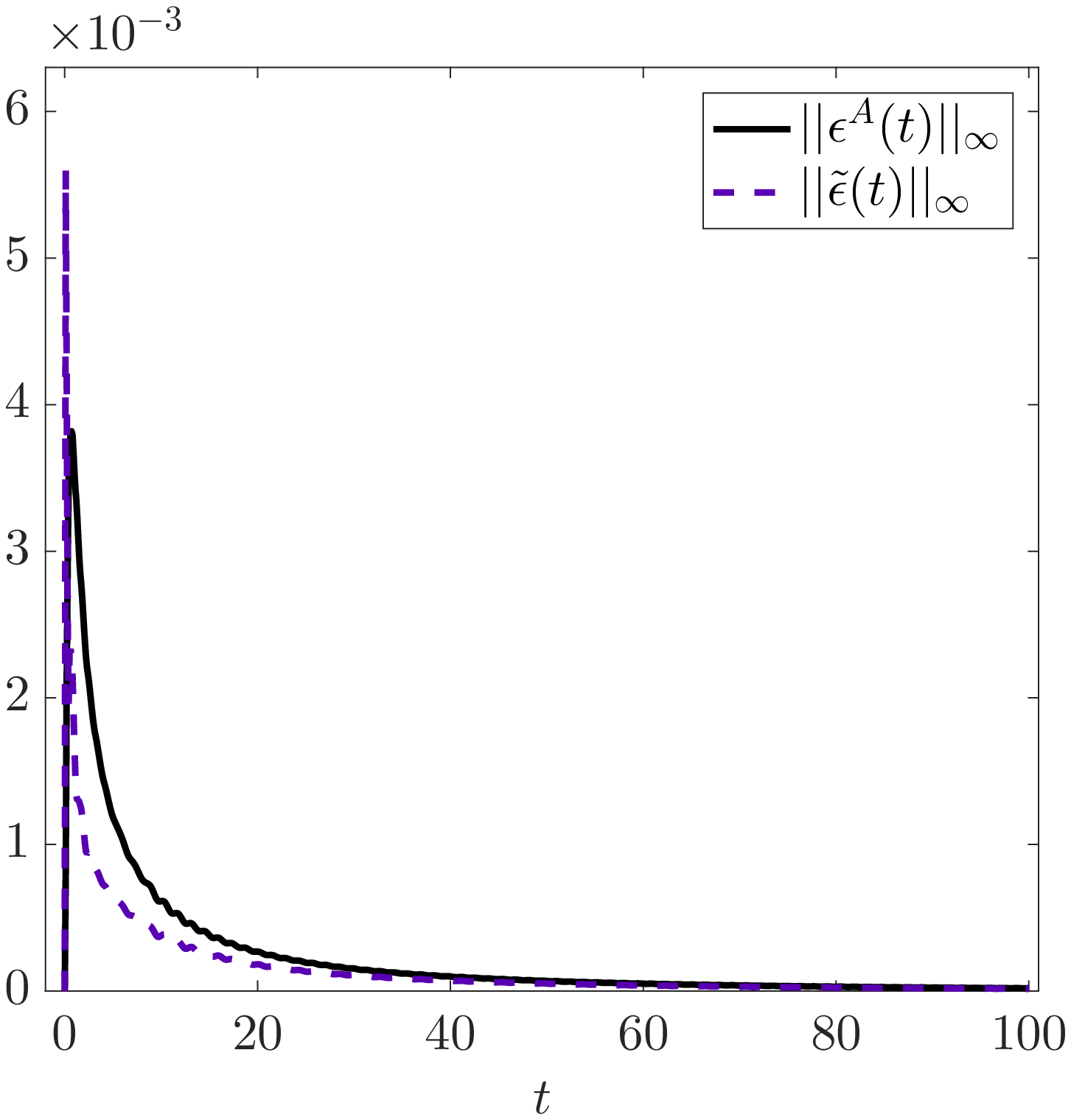}
\label{fig:QueueLengthDistModStepUpLongTimeErr}
}
\caption{Comparison of queue-length distribution with exact ODE system and continuous approximation in the moderate ramp up case, (a) probability $p_0(t)$, (b) $p_3(t)$, (c) supremum error.}
\label{fig:QueueLengthDistModStepUpLongTime}
\end{figure}

Generally, we have two types of steps, the ramp up and the ramp down, where we interpret up and down by the value of the traffic intensity $\varrho$. We first consider three types of ramp ups: moderate, strong and very strong ramp ups. A moderate ramp up is determined by $\lambda_0 = 0.5$, $\lambda_1 = 0.8$ and $\mu_0 = \mu_1 = 1$, which corresponds to an increase in the traffic intensity $\varrho_0 = 0.5$ to $\varrho_1 = 0.8$. In figure \ref{fig:QueueLengthDistModStepUpLongTime} (a) - (b), the values of $p_k(t)$, $p_k^A(t)$ and $\tilde{p}_k(t)$ are shown for different time points and $k \in \{0,3\}$. Visually, the continuous approximations given by the squares and diamond markers are coincident with the exact model given by the black dots. Some small displacements can  be observed at the second and third time points. 
This is emphasized by the numerical error measures displayed in figure \ref{fig:QueueLengthDistModStepUpLongTimeErr}. As expected from the derivation of the continuous approximation, we have a decay of the error as time evolves since the model is exact in the steady state. The largest errors occur right after the step at $t=0$, which is intuitive since it is the time right at the disturbance  and also corresponds to the observation that the analytic formulas \eqref{eq:rhoStSt1}  and \eqref{eq:PStSt1} are evaluated at singularities for $t \to 0$, which implies errors. Nevertheless, the maximal difference in the exact model is  small,  of order $10^{-3}$, as seen in figure \ref{fig:QueueLengthDistModStepUpLongTime} (c). We relate ``small'' to the values given in figure \ref{fig:QueueLengthDistModStepUpLongTime} (a)-(b).

\begin{figure}[htb!]
%
\subfigure[Approximation of $p_0(t)$]{
\includegraphics[width=0.3\textwidth]{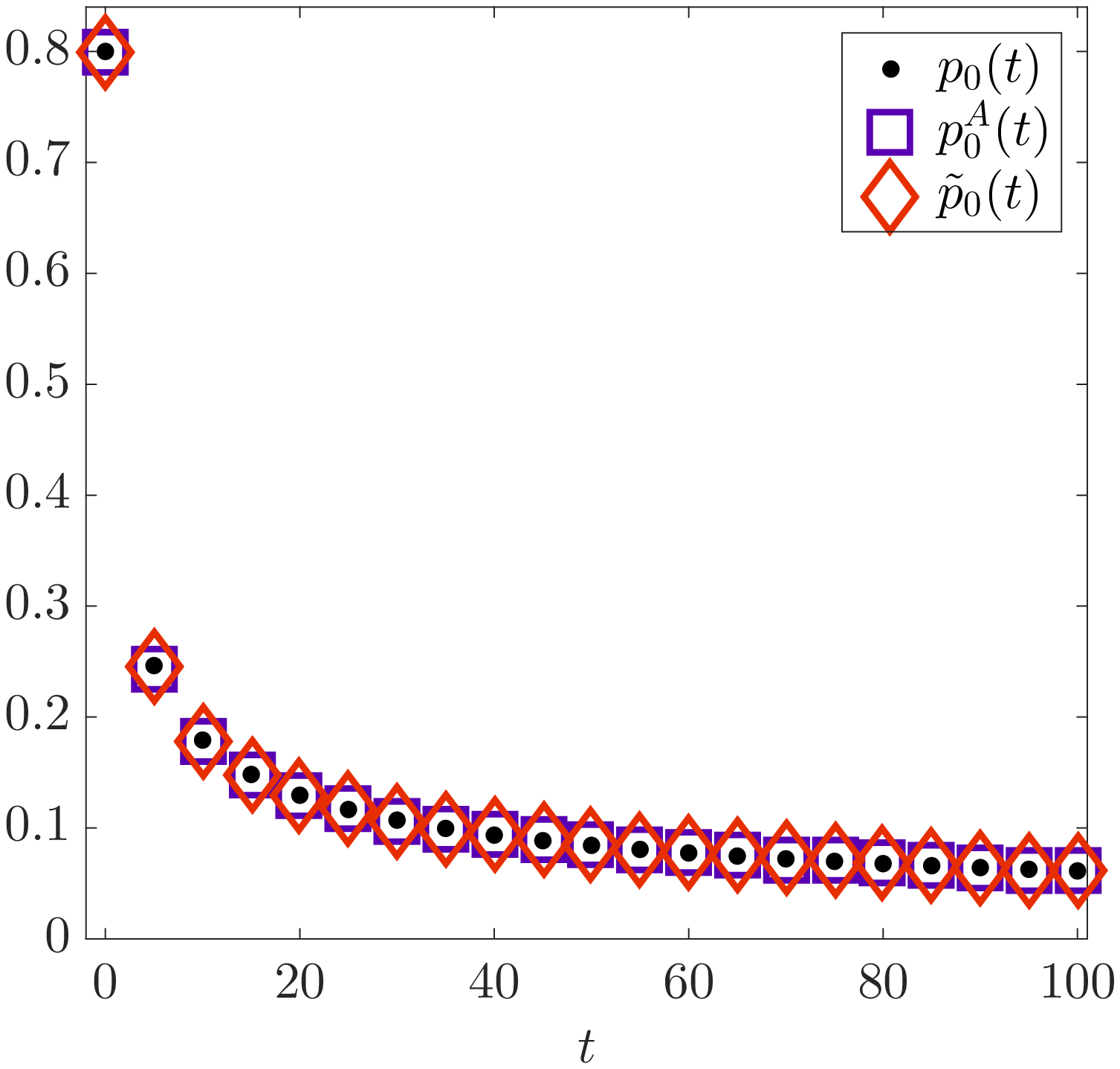}
}
\subfigure[Approximation of $p_3(t)$]{
\includegraphics[width=0.3\textwidth]{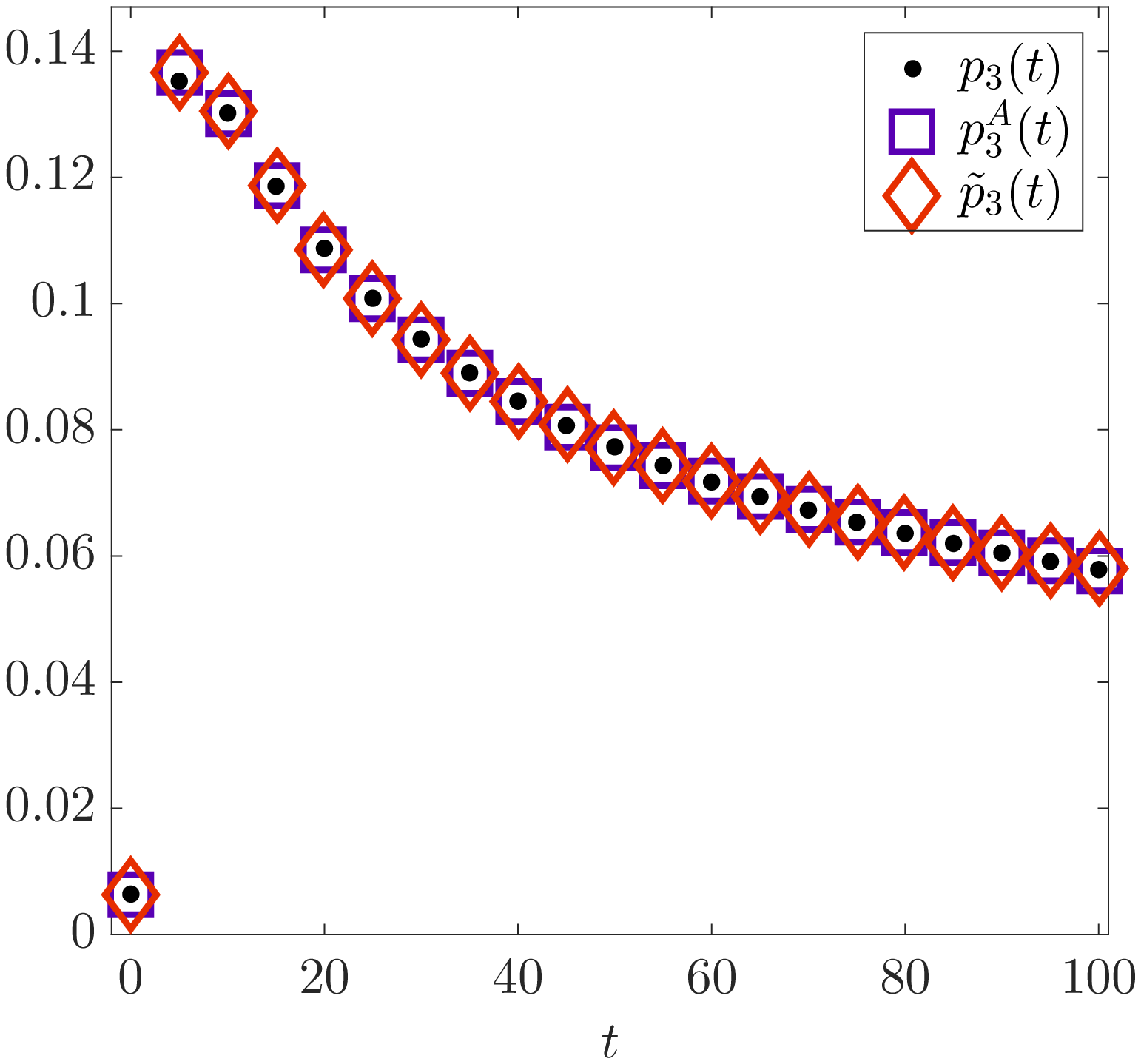}
}
\subfigure[Supremum error]{
\includegraphics[width=0.3\textwidth]{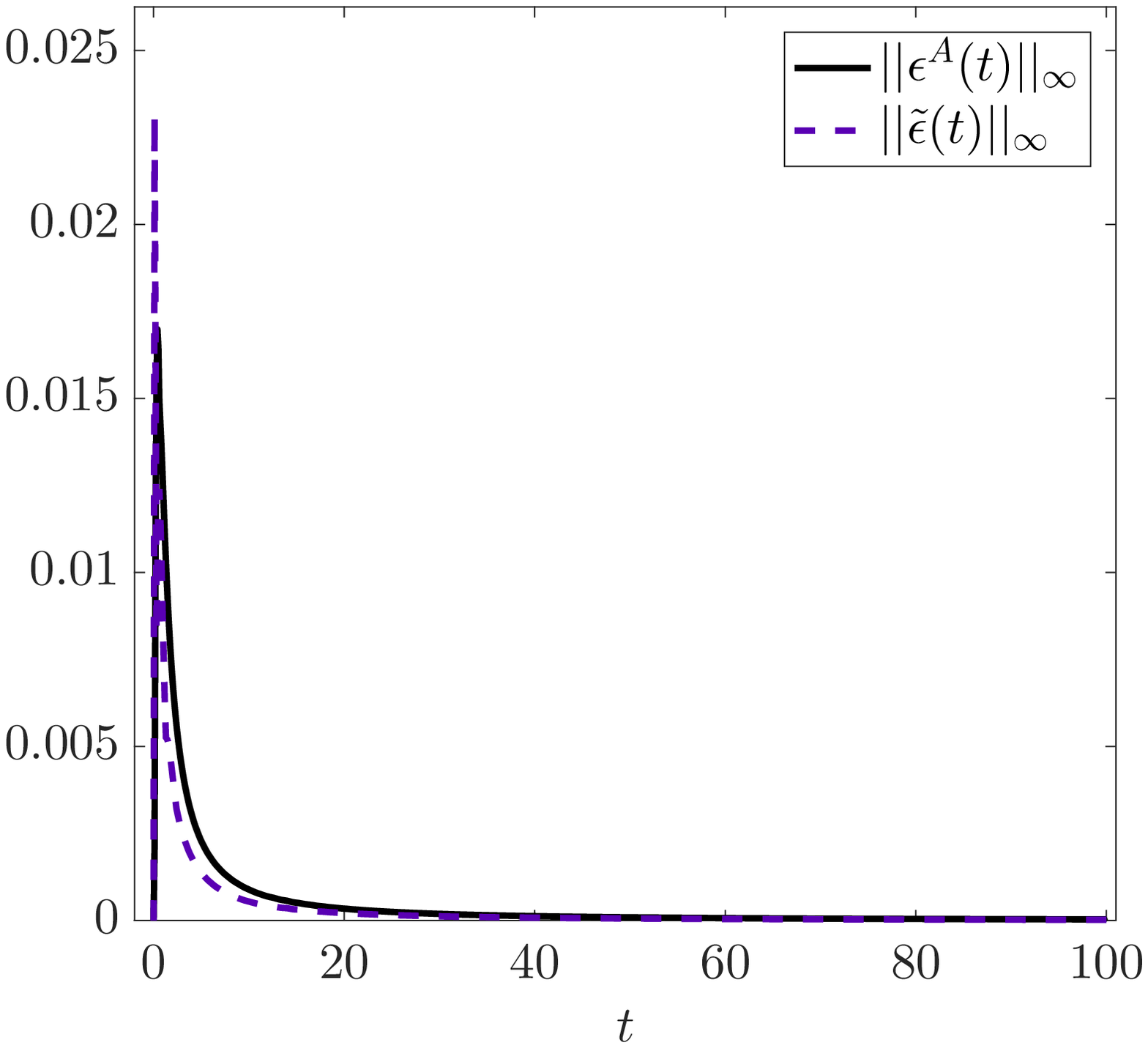}
\label{fig:QueueLengthDistStrongStepUpLongTimeErr}
}
\caption{Comparison of queue-length distribution with exact ODE system and continuous approximation in the strong ramp up case, (a) probability $p_0(t)$, (b) $p_3(t)$, (c) supremum error.}

\label{fig:QueueLengthDistStrongStepUpLongTime}
\end{figure}

Considering  a strong ramp up given by $\lambda_0 = 0.2$, $\lambda_1 = 0.99$ and $\mu_0 = \mu_1 = 1$, we again obtain a useful approximation, as shown in figure \ref{fig:QueueLengthDistStrongStepUpLongTime} (a)-(b). The displacements are slightly larger in the first time steps, and the errors again decrease in time; see figure \ref{fig:QueueLengthDistStrongStepUpLongTime} (c). Since $\varrho_1 = 0.99<1$, the system still converges to a steady state. The supremum error is of order $10^{-2}$ in this case, where the largest errors occur again right after the step; they are small compared to the values in figure \ref{fig:QueueLengthDistStrongStepUpLongTime} (a)-(b).

\begin{figure}[H]
\subfigure[Approximation of $p_0(t)$]{
\includegraphics[width=0.35\textwidth]{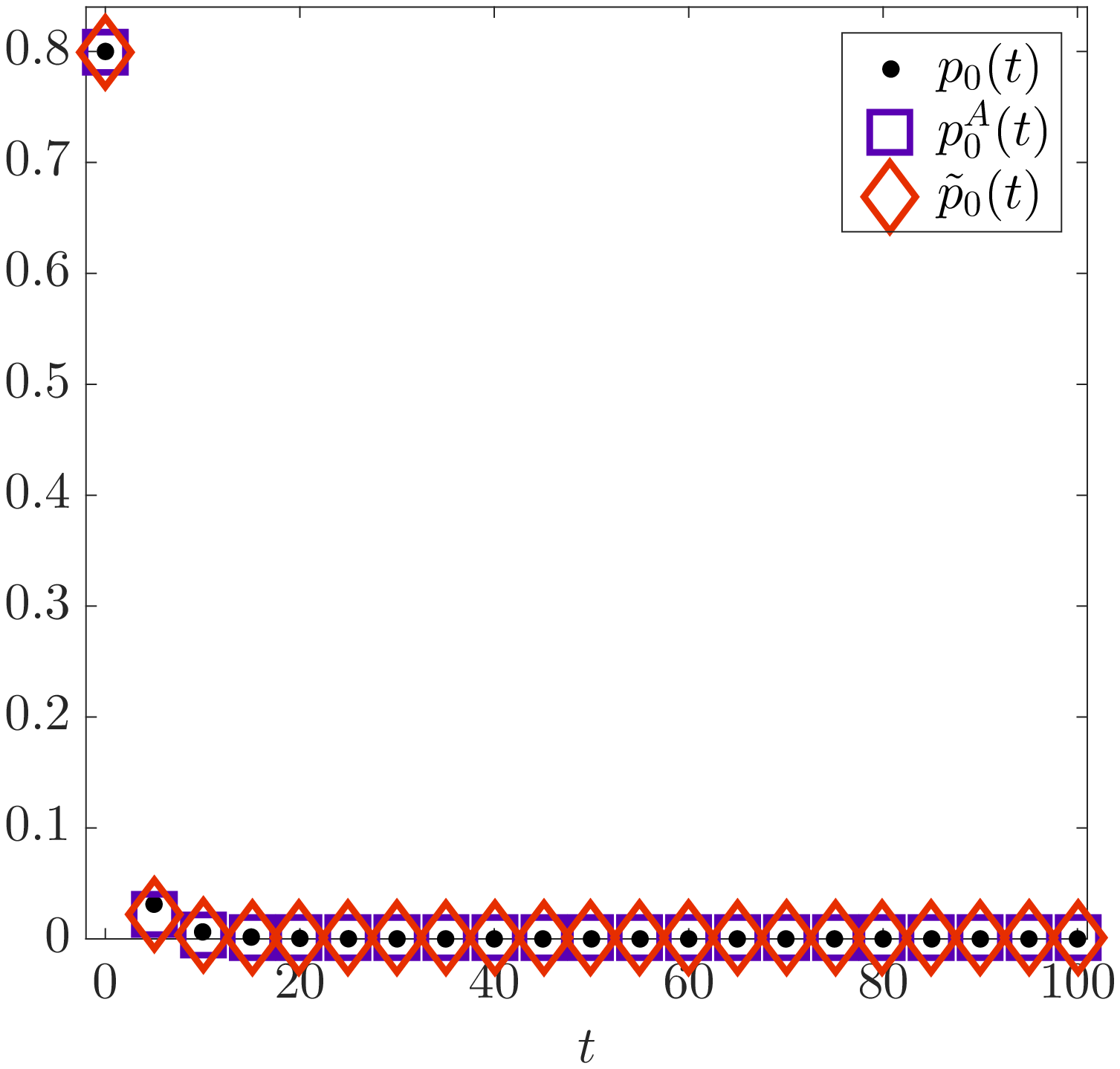}
}
\subfigure[Approximation of $p_{10}(t)$]{
\includegraphics[width=0.35\textwidth]{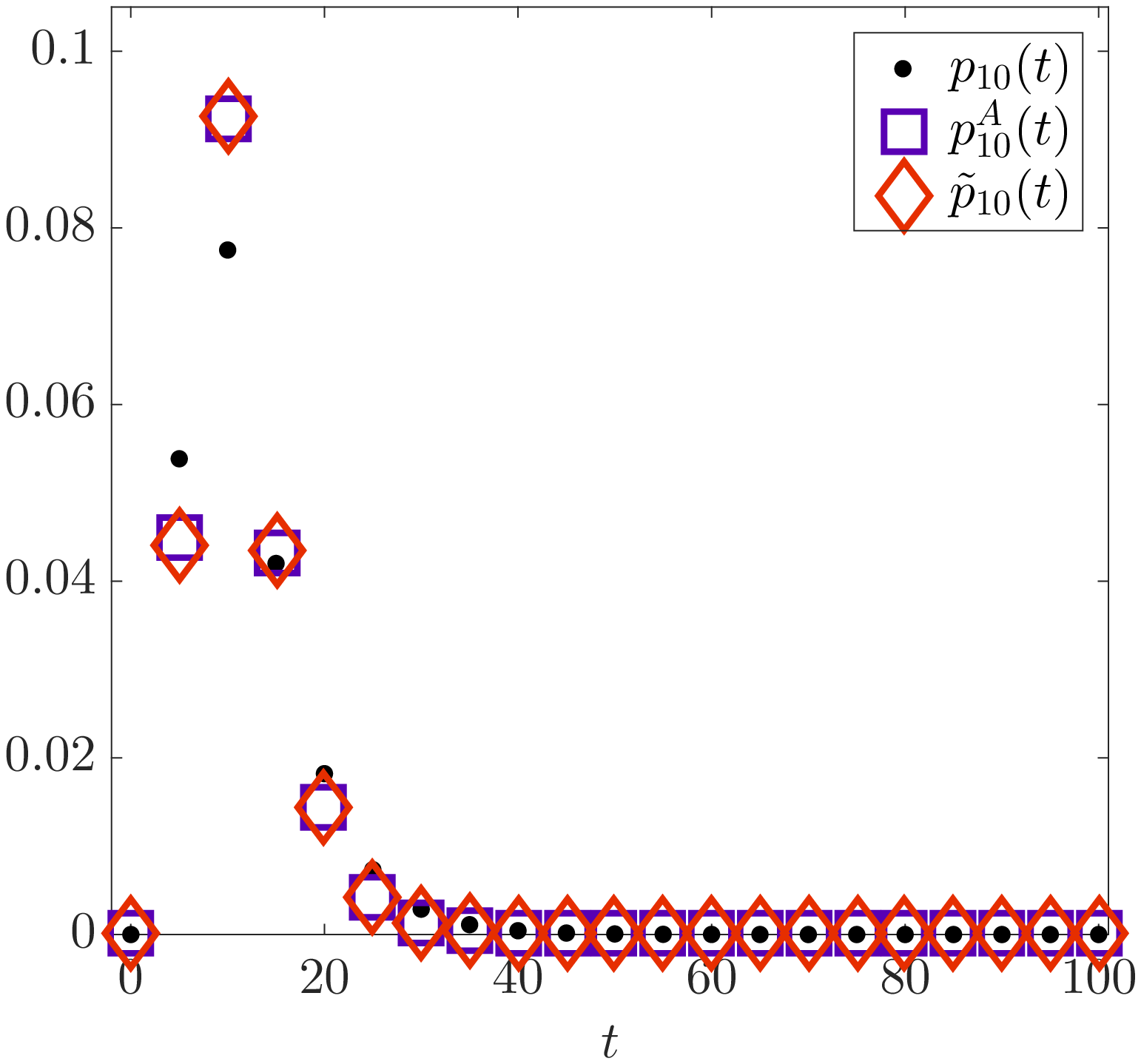}
}
\\
\subfigure[Approximation of $p_{40}(t)$]{
\includegraphics[width=0.35\textwidth]{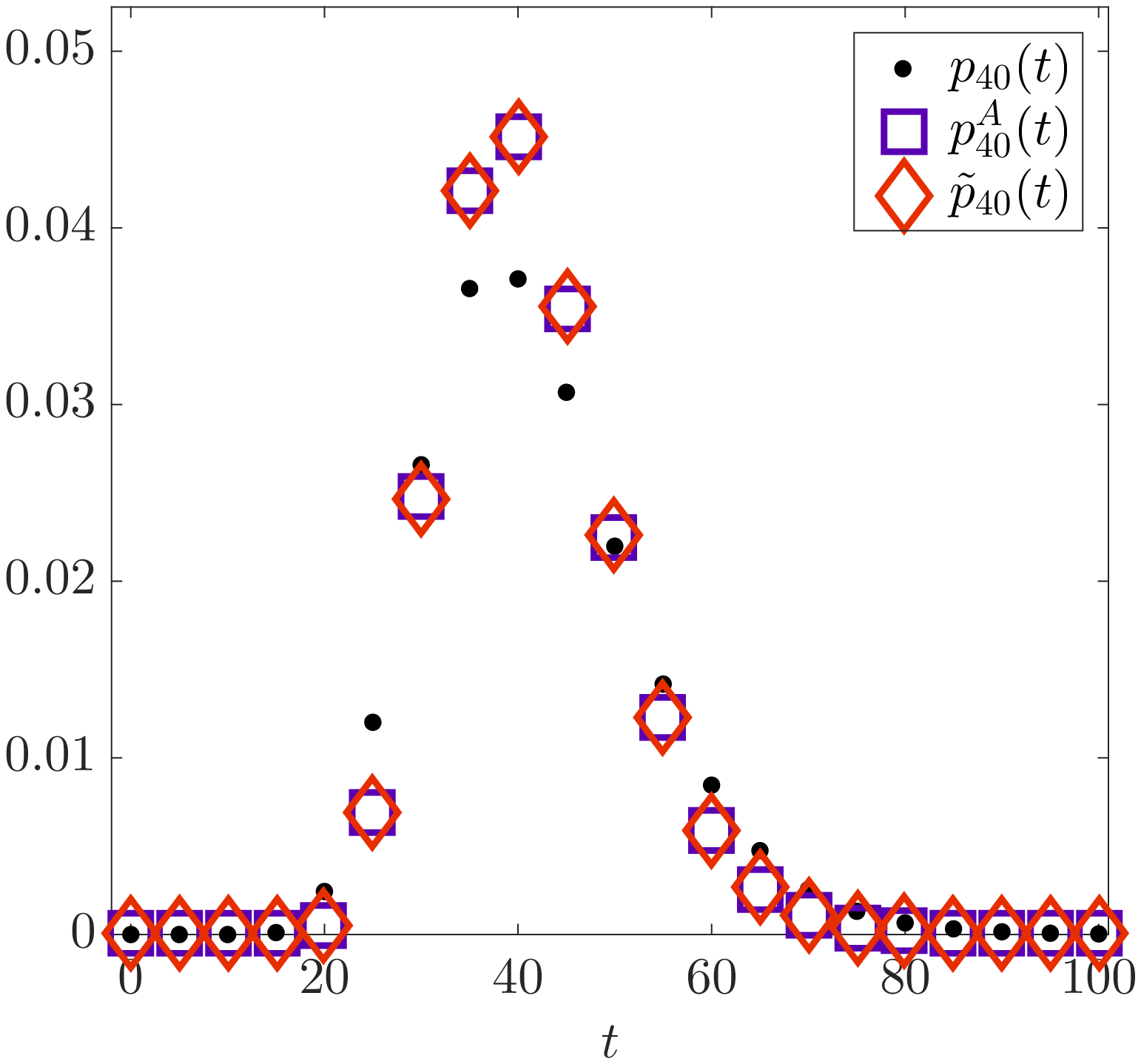}
}
\subfigure[Supremum error]{
\includegraphics[width = 0.35\textwidth]{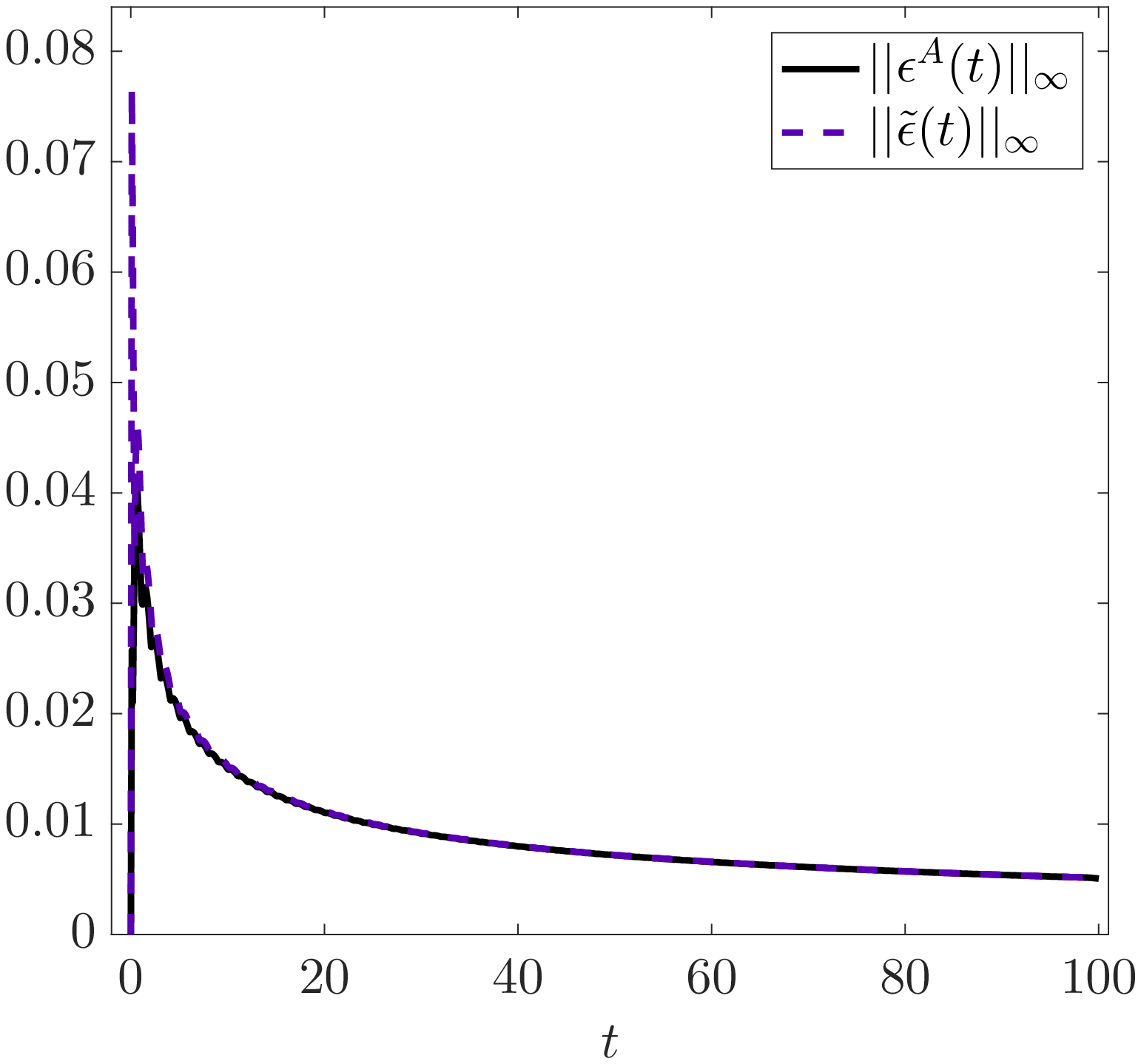}
\label{fig:QueueLengthDistVeryStrongStepUpLongTimeErr}
}
\caption{Comparison of queue-length distribution with exact ODE system and continuous approximation in the very strong ramp up case, 
(a) probability $p_0(t)$, (b) $p_{10}(t)$, (c) $p_{40}(t)$, (d) supremum error.}
\label{fig:QueueLengthDistVeryStrongStepUpLongTime}
\end{figure}

The last ramp up example that we use is a very strong ramp up from $\lambda_0 = 0.2$ to $\lambda_1 = 2$  with $\mu_0 = \mu_1 = 1$. In this case, the system has no steady state, which is a fundamental assumption in the derivation of the continuous approximation. We again observe a decreasing supremum error in figure \ref{fig:QueueLengthDistVeryStrongStepUpLongTime} (d); however, compared to the values given in \ref{fig:QueueLengthDistVeryStrongStepUpLongTime} (a)-(c), we cannot deduce that they are small. 
Of course, a constant utilization in time greater than one is not a realistic scenario.
We analyze the approximation for short times of over utilization ($\varrho(t)\geq 1$) in the time-dependent case later.

\subsubsection*{Ramp Down}
\begin{figure}[H]
\subfigure[Approximation of $p_0(t)$]{
\includegraphics[width = 0.3\textwidth]{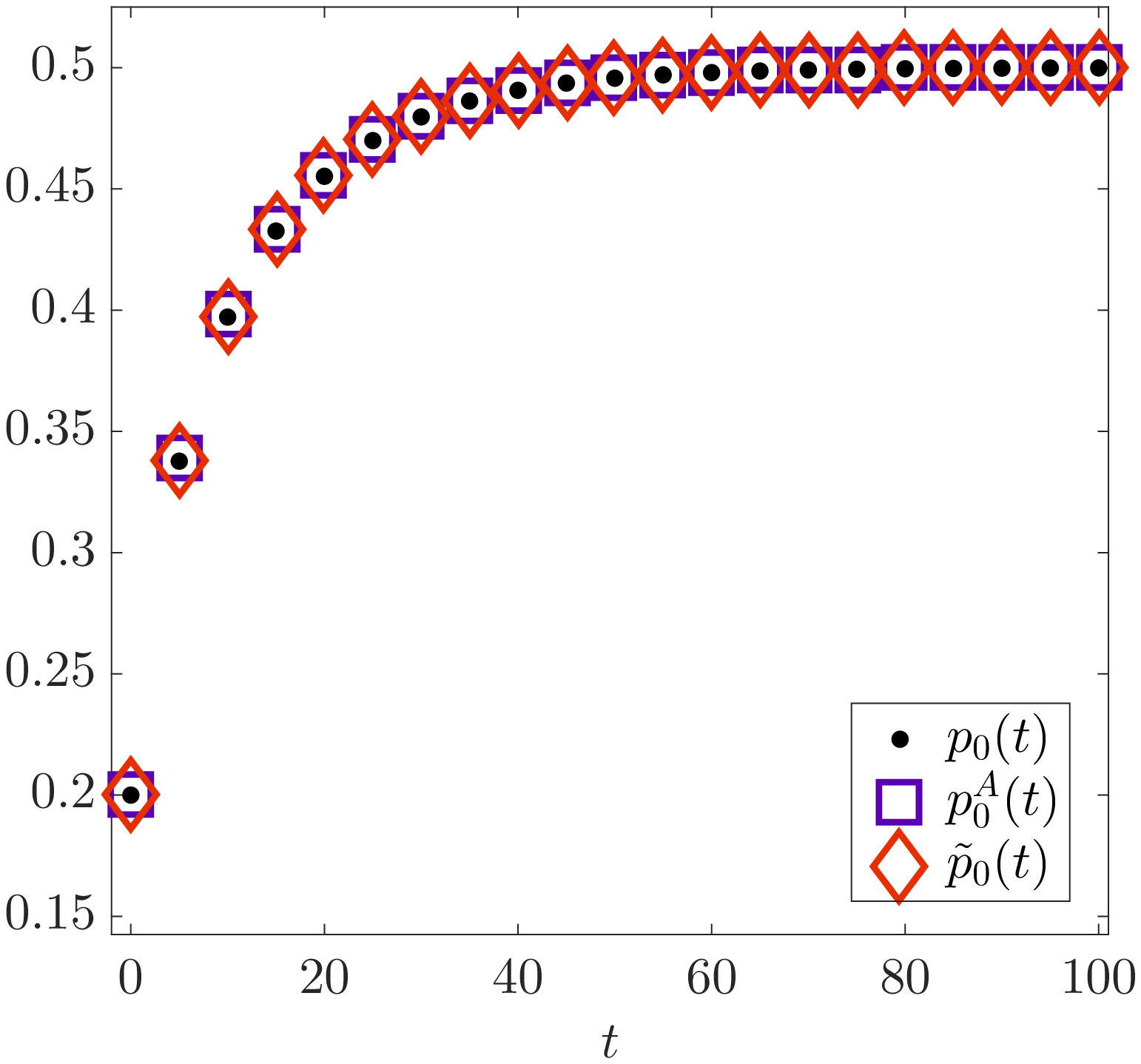}
}
\subfigure[Approximation of $p_3(t)$]{
\includegraphics[width = 0.3\textwidth]{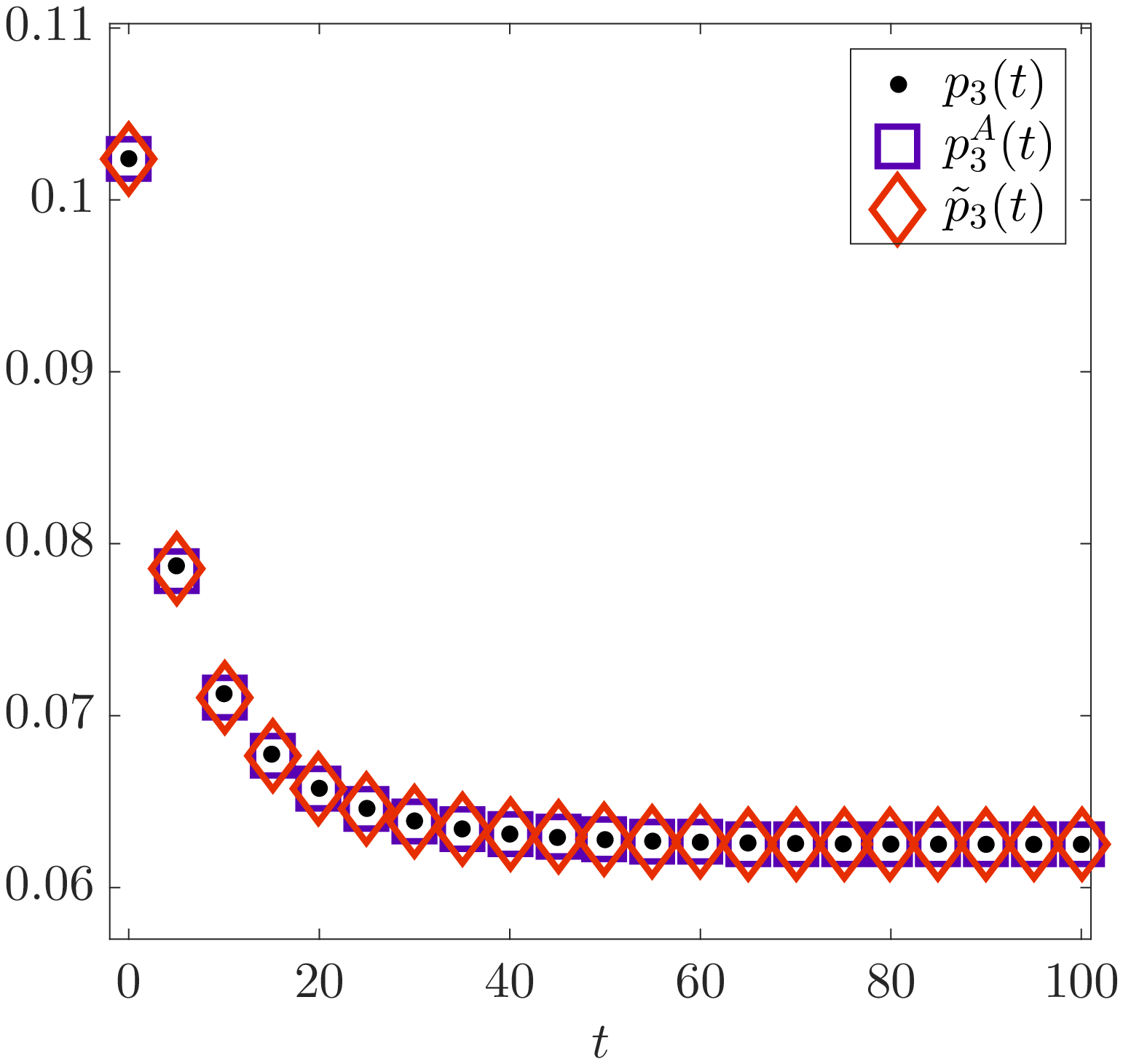}
}
\subfigure[Supremum error]{
\includegraphics[width = 0.3\textwidth]{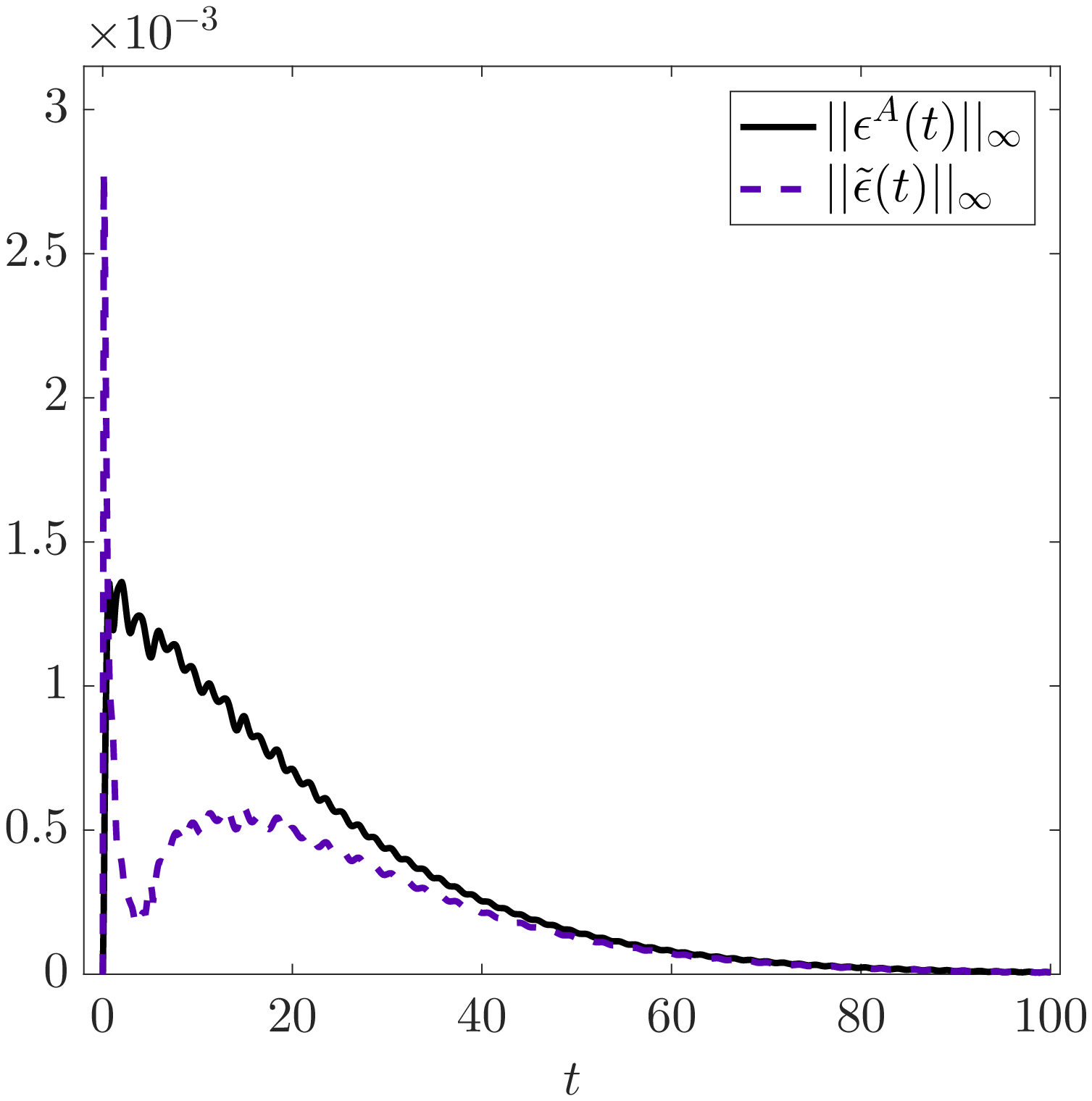}
\label{fig:QueueLengthDistModStepDownLongTimeErr}
}
\caption{Comparison of queue-length distribution with exact ODE system and continuous approximation in the moderate ramp down case, (a) probability $p_0(t)$, (b) $p_3(t)$, (c) supremum error.}

\label{fig:QueueLengthDistModStepDownLongTime}
\end{figure}

Since we assume a steady state at $t=0$, we only consider a moderate and strong ramp down in the following. The moderate ramp down is given by $\lambda_0 = 0.8$, $\lambda_1 = 0.5$ and $\mu_0 = \mu_1 = 1$.  Visually, the approximations are exact as  seen in figure \ref{fig:QueueLengthDistModStepDownLongTime} (a)-(b), and  supported by the error measures in figure \ref{fig:QueueLengthDistModStepDownLongTime} (c), which are of order $10^{-3}$. 
%
We see a decreasing error in time again; only the approximation $\tilde{p}_k(t)$ shows a different behavior initially. One reason for this phenomenon is the evaluations of the standard normal cdf at singular values, which are different in the case of $\tilde{p}_k(t)$.

\begin{figure}[h]
\subfigure[Approximation of $p_0(t)$]{
\includegraphics[width = 0.3\textwidth]{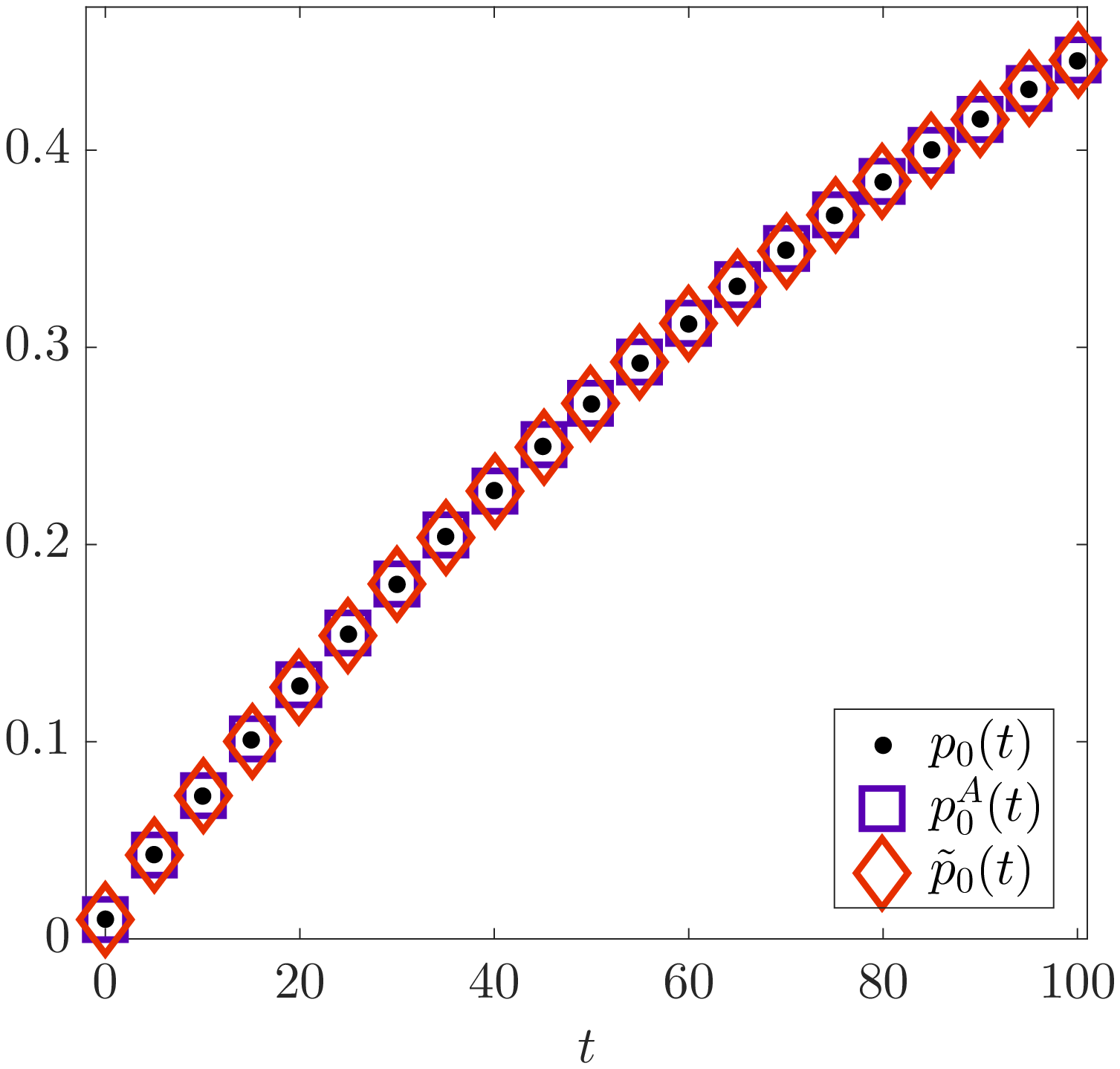}
}
\subfigure[Approximation of $p_3(t)$]{
\includegraphics[width = 0.3\textwidth]{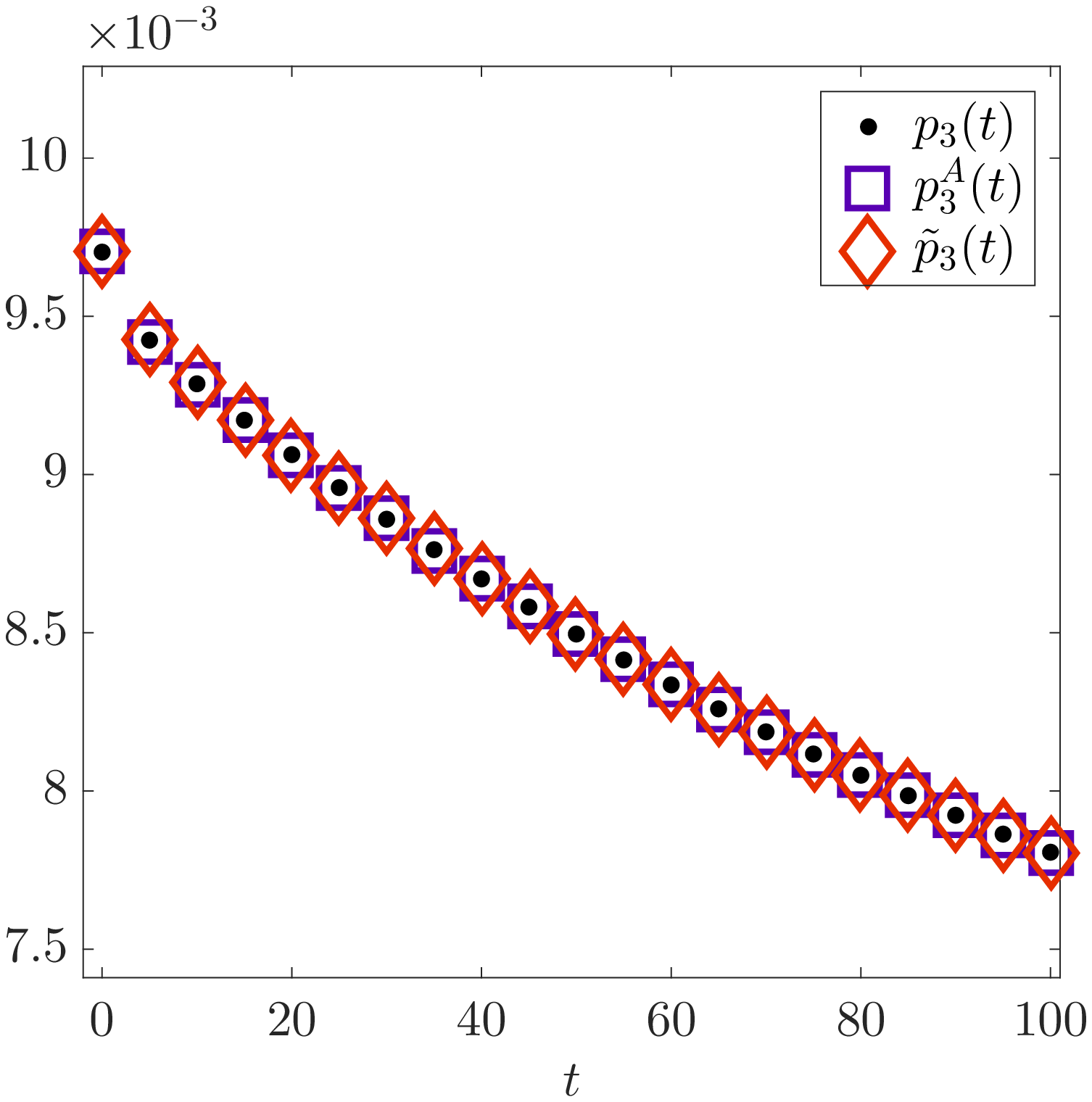}
}
\subfigure[Supremum error]{
\includegraphics[width = 0.3\textwidth]{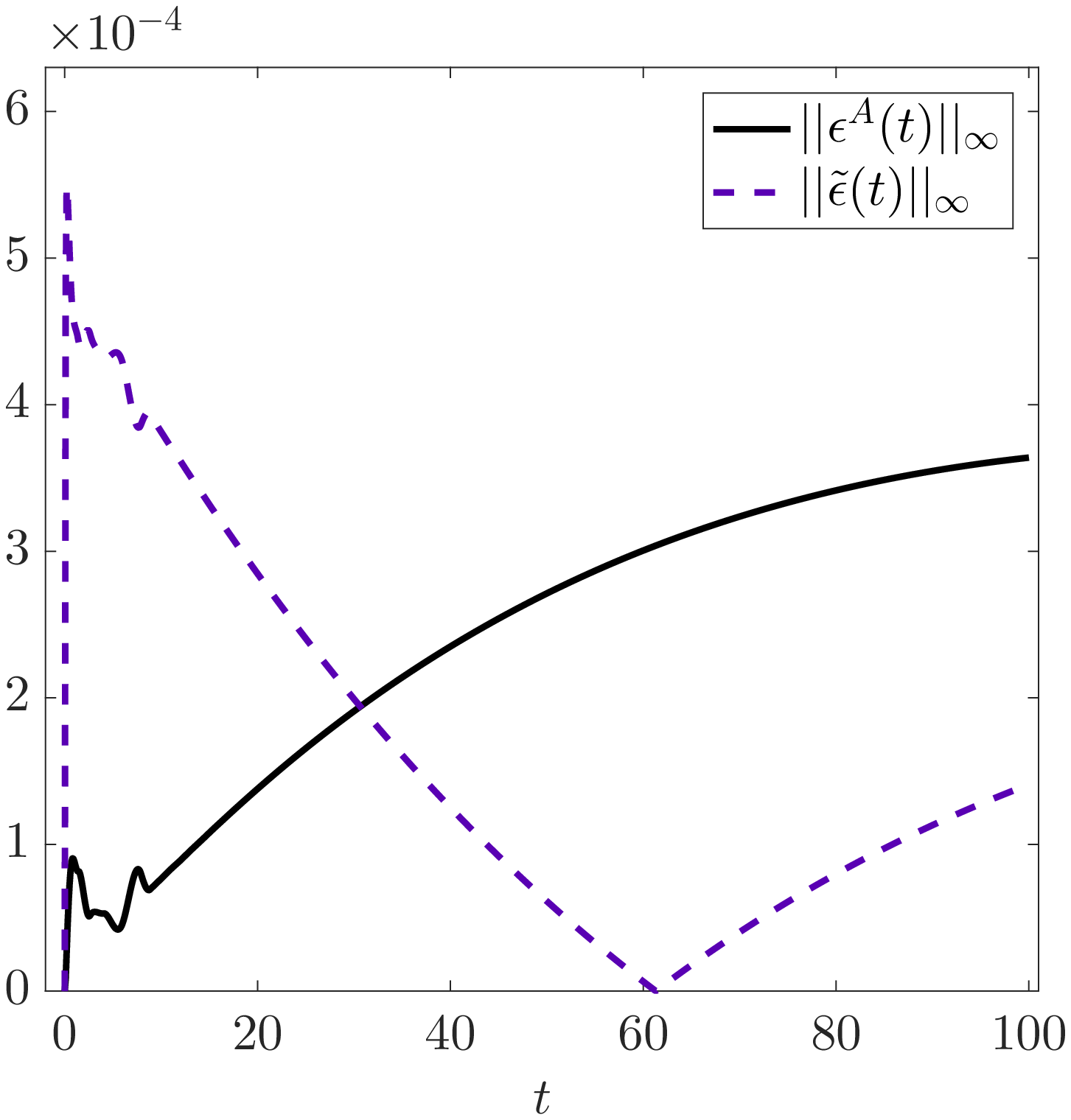}
\label{fig:QueueLengthDistStrongStepDownLongTimeErr}
}
\caption{Comparison of queue-length distribution with exact ODE system and continuous approximation in the strong ramp down case, (a) probability $p_0(t)$, (b) $p_3(t)$, (c) supremum error.}
\label{fig:QueueLengthDistStrongStepDownLongTime}
\end{figure}

In the case of a strong ramp down $\lambda_0 = 0.99$, $\lambda_1 = 0.2$ and $\mu_0 = \mu_1 = 1$, we have a quite good approximation, as shown in figure \ref{fig:QueueLengthDistStrongStepDownLongTimeErr} (a)-(b). The errors in figures \ref{fig:QueueLengthDistStrongStepDownLongTimeErr} (c) are of order $10^{-4}$, and visually, the values $p_k(t)$ for $k \in\{ 0, \dots 3\}$ are fairly well approximated.

%
%
%

\subsubsection*{Time-dependent Coefficient: Cyclic}
We consider a cyclic time-varying inflow rate, i.e., 
\begin{align*}
\lambda(t) = \frac{\gl_0-\gl_1}{2}\cos\left(2 \pi \frac{t}{T_{\text{Per}}}\right)+\frac{\lambda_0+\gl_1}{2}
\end{align*}
 as has been studied in \cite{Wienke2015} to approximate the expected outflow of a queueing system with a second-order model of hyperbolic equations. The parameter $\gl_0\geq 0$ denotes the lowest and $\gl_1\geq 0$ the highest value of the inflow rate, which is periodic with period $T_{\text{Per}}>0$. We set the production rate as a constant of $\mu(t) \equiv 1$ and, analogously to the step case, study different values for $\gl_0$ and $\gl_1$.
The probabilities $p_k(t)$ are again computed with the ODE system, and the approximate values $p_k^A(t)$ are approximated with the numerical scheme \eqref{eq:NumSchemeAppQueue1}-\eqref{eq:NumSchemeAppQueue2}, where we use $\Delta x = 0.02$, $x_0 = 0$, $x_1 = 200$ and a time horizon $T=25$. The integration of $p_k^A(t)$ is done with a trapezoidal rule.

\begin{figure}[H]
\subfigure[Approximation of $p_0(t)$, $p_1(t)$ for $T_{\text{Per}} = 10$]{
\includegraphics[width = 0.4\textwidth]{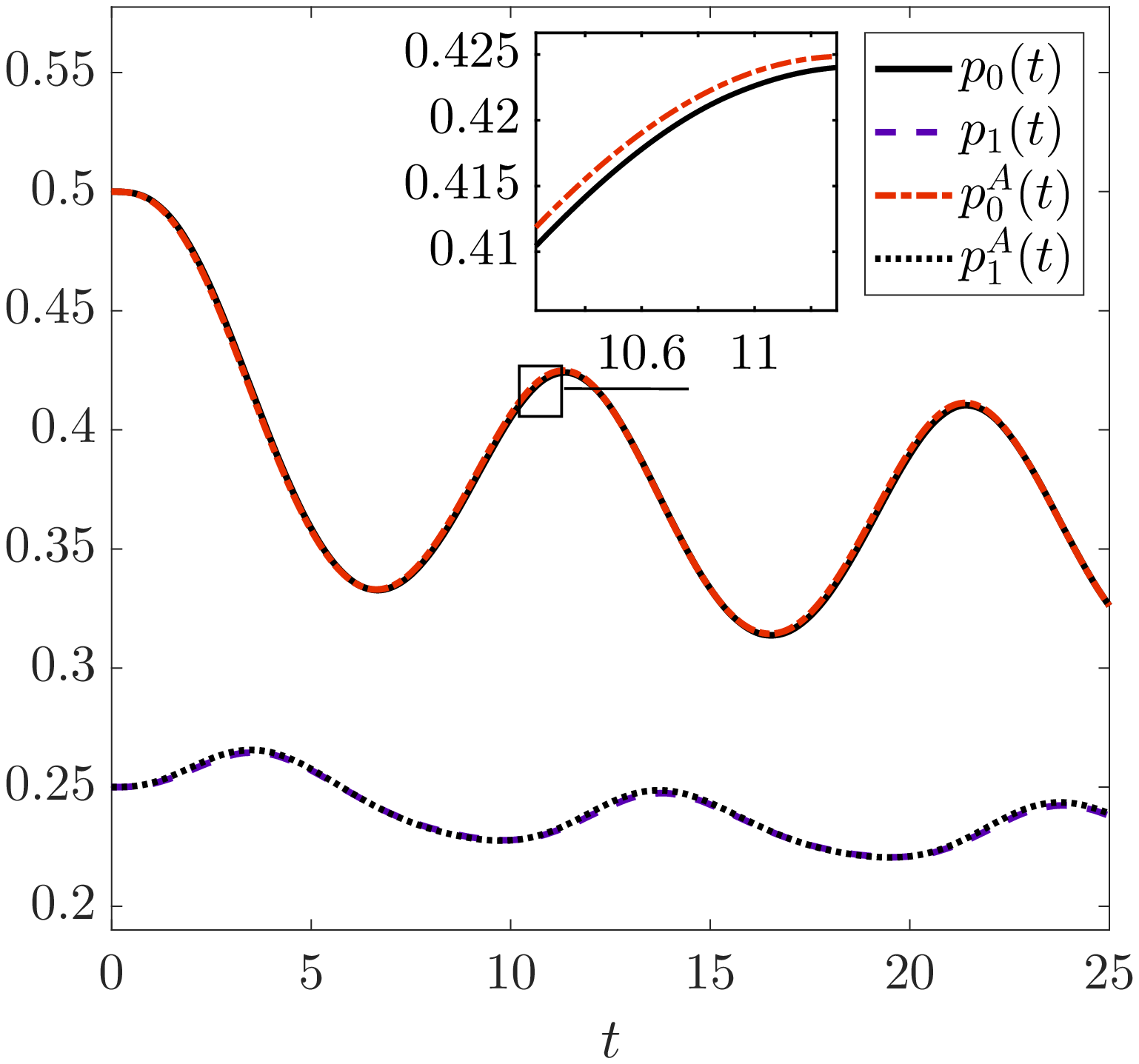}
}
\subfigure[Approximation of $p_0(t)$, $p_1(t)$ for $T_{\text{Per}} = 2$]{
\includegraphics[width = 0.4\textwidth]{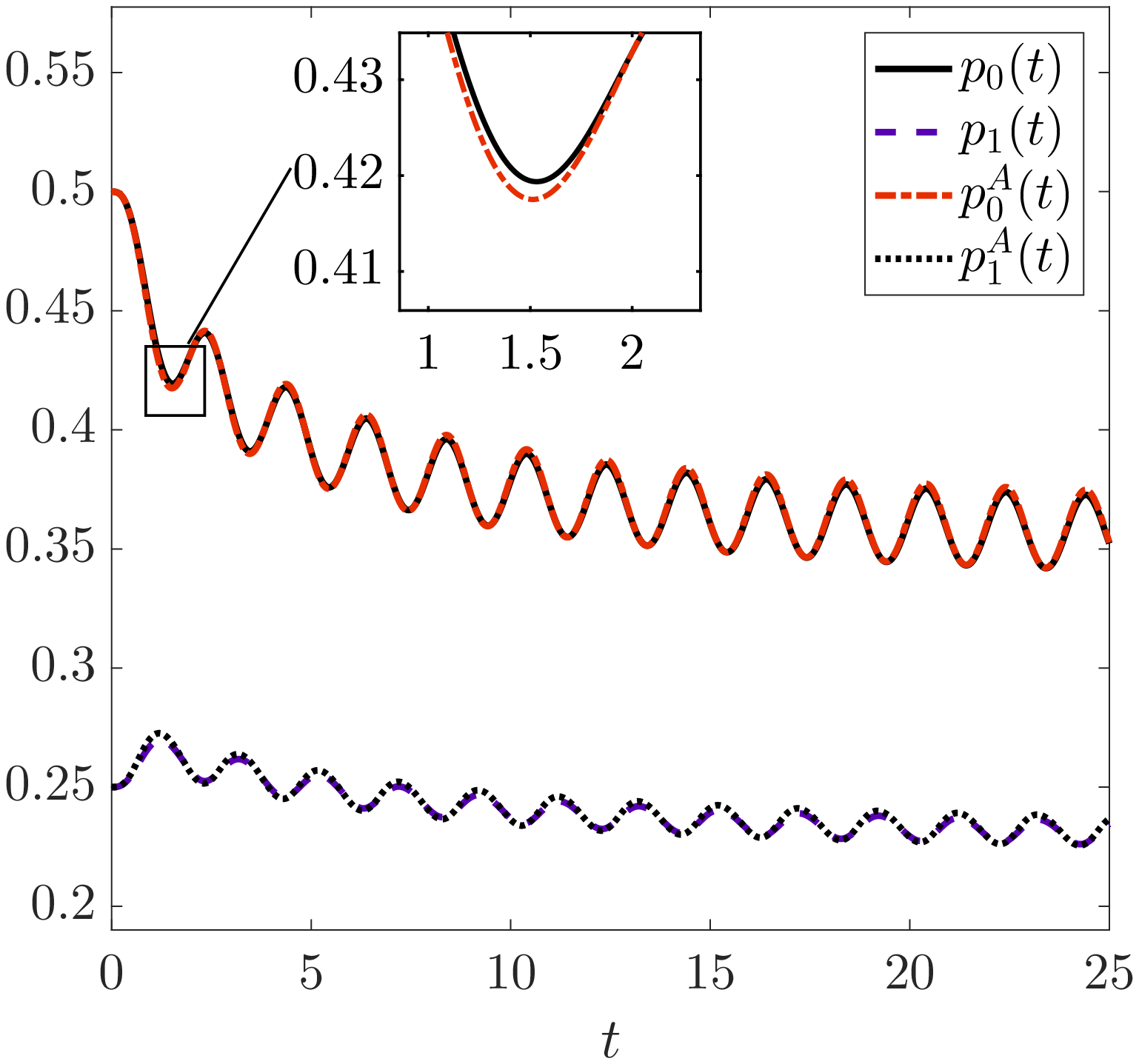}
}
\caption{Comparison of queue-length distribution with exact ODE system and continuous approximation in the moderate cyclic inflow rate case}
\label{fig:QueueLengthDistModerateCyclic}
\end{figure}

We again study a moderate, strong and very strong case and  call the case $\gl_0 = 0.5$ and $\gl_1 = 0.8$ the moderate case. In figure \ref{fig:QueueLengthDistModerateCyclic}, we show the results of a simulation for the probability that no and one customer are in the queueing system, respectively. For both periods,  $T_{\text{Per}} = 10$ and $T_{\text{Per}} = 2$, visually, the approximations $p_0^A(t)$ and $p_1^A(t)$ are close to the values $p_0(t)$ and $p_1(t)$.
Table \ref{tab:ErrorTimeDepQueueingModerateCyclic} first column shows the maximal difference between the approximation and the ODE result for $k = 0,\dots,100$. We observe that the supremum norm of the error increases with a smaller period $T_{\text{Per}}$ but remains of the order $10^{-3}$, which is small compared to the values in figure \ref{fig:QueueLengthDistModerateCyclic}.

\begin{table}[H]
\centering
\begin{tabular}{c||c|c|c}
$T_{\text{Per}}$ & &$\displaystyle \max_j \|\epsilon^A(t_j)\|_\infty \;[10^{-3}]$  \\ 
\hline 
& moderate cyclic case & strong cyclic case & very strong cyclic case \\
\hline \hline
25 & 1.4876 & 6.8495  & 6.5022 \\ 
\hline 
10 & 1.7911 & 8.8031 & 11.6854\\ 
\hline 
5 & 2.2439 & 12.0694 & 17.2880 \\ 
\hline 
2 & 3.1356 & 17.1375 & 25.8903 \\ 
\hline 
1 & 3.2942 &18.1220 & 30.9673  \\ 
\end{tabular}
\caption{Error of the continuous approximation for different periods in the moderate, strong and very strong  cyclic case}
\label{tab:ErrorTimeDepQueueingModerateCyclic} 
\end{table}

\begin{figure}[hbt!]
\subfigure[Approximation of $p_0(t)$, $p_1(t)$ for $T_{\text{Per}} = 10$]{
\includegraphics[width = 0.4\textwidth]{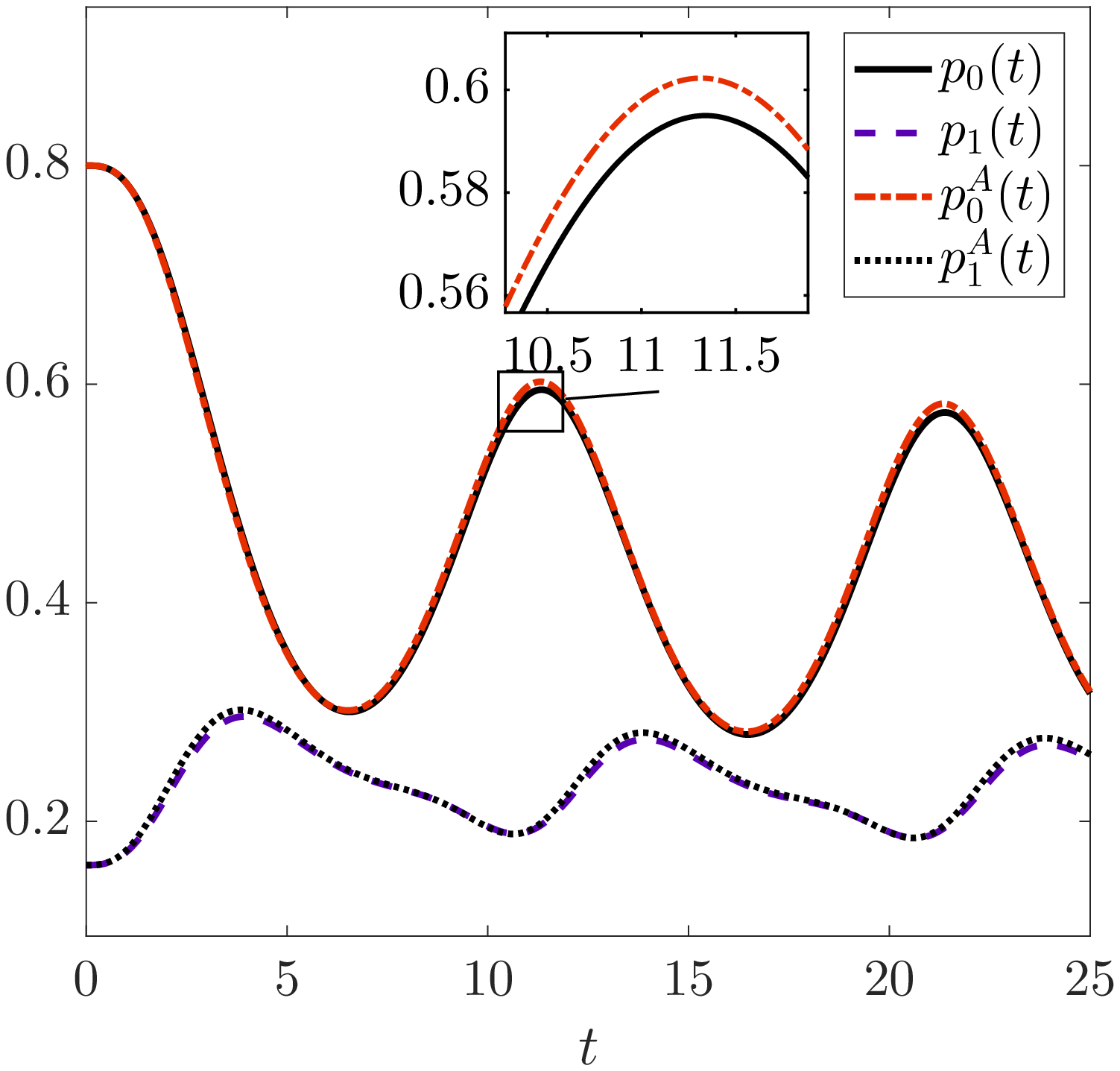}
}
\subfigure[Approximation of $p_0(t)$, $p_1(t)$ for $T_{\text{Per}} = 2$]{
\includegraphics[width = 0.4\textwidth]{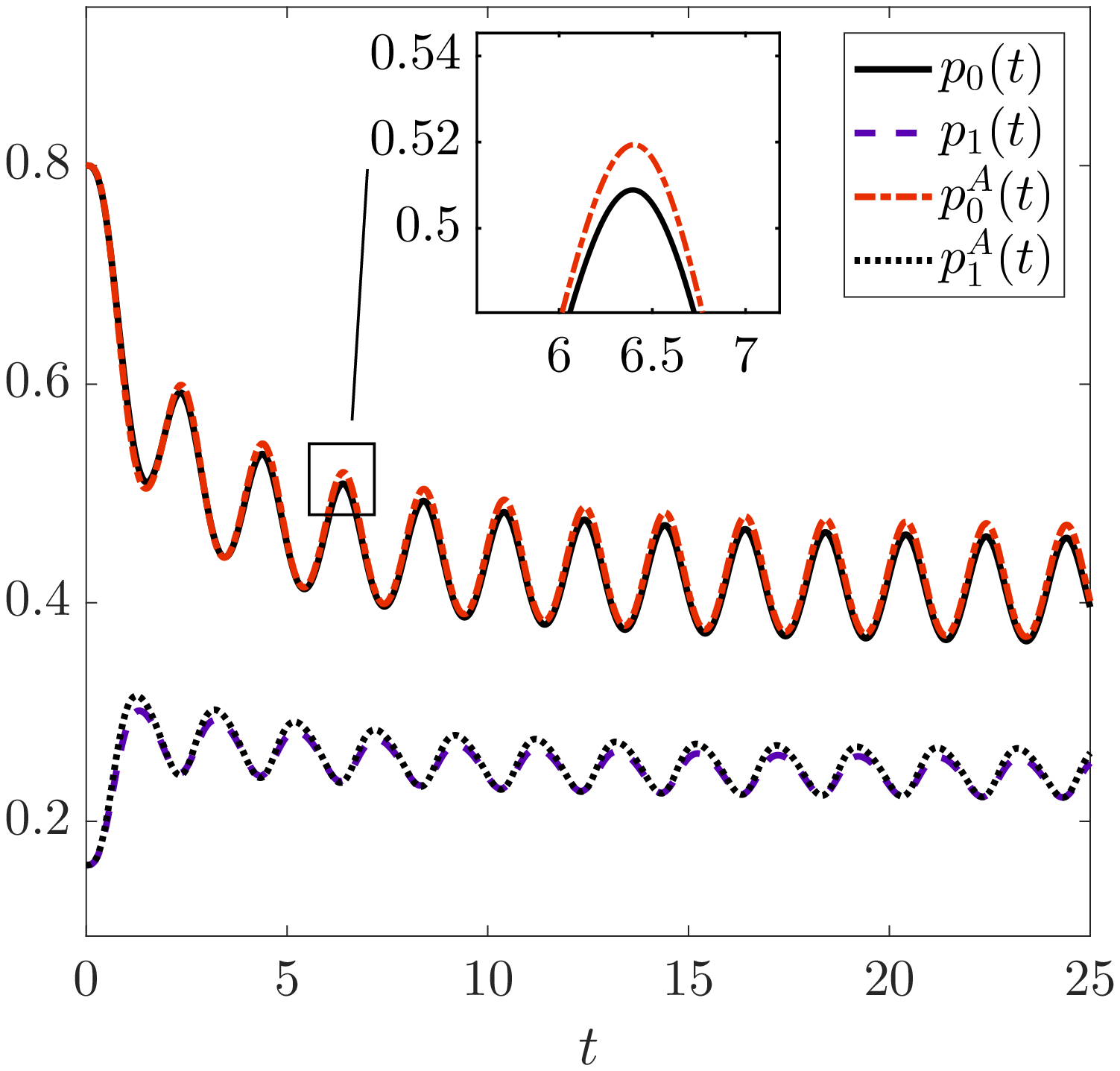}
}
\caption{Comparison of queue-length distribution with exact ODE system and continuous approximation in the strong cyclic inflow rate case}
\label{fig:QueueLengthDistStrongCyclic}
\end{figure}

In the strong case given by $\gl_0 = 0.2$ and $\gl_1 = 0.99$ corresponding to larger amplitude oscillations, we can find in  figure \ref{fig:QueueLengthDistStrongCyclic} a larger deviation of the approximate model from the ODE system than in the moderate case. Table \ref{tab:ErrorTimeDepQueueingModerateCyclic} second column shows that the numerical error measures  are of order $10^{-2}$ in this case and are again increasing as $T_{\text{Per}}$ decreases.

\begin{figure}[hbt!]
\subfigure[Approximation of $p_0(t)$, $p_1(t)$ for $T_{\text{Per}} = 10$]{
\includegraphics[width = 0.4\textwidth]{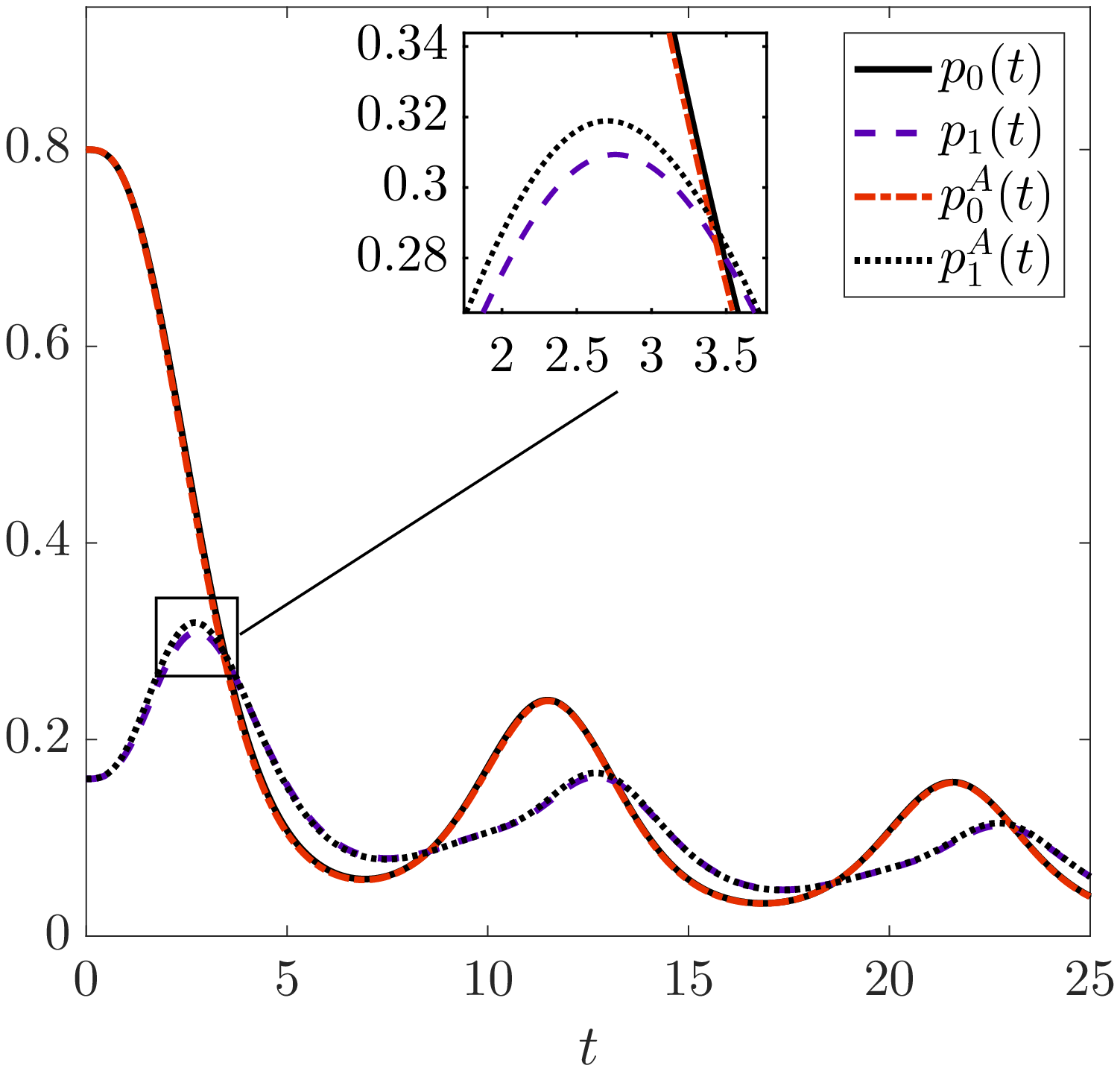}
}
\subfigure[Approximation of $p_0(t)$, $p_1(t)$ for $T_{\text{Per}} = 2$]{
\includegraphics[width = 0.4\textwidth]{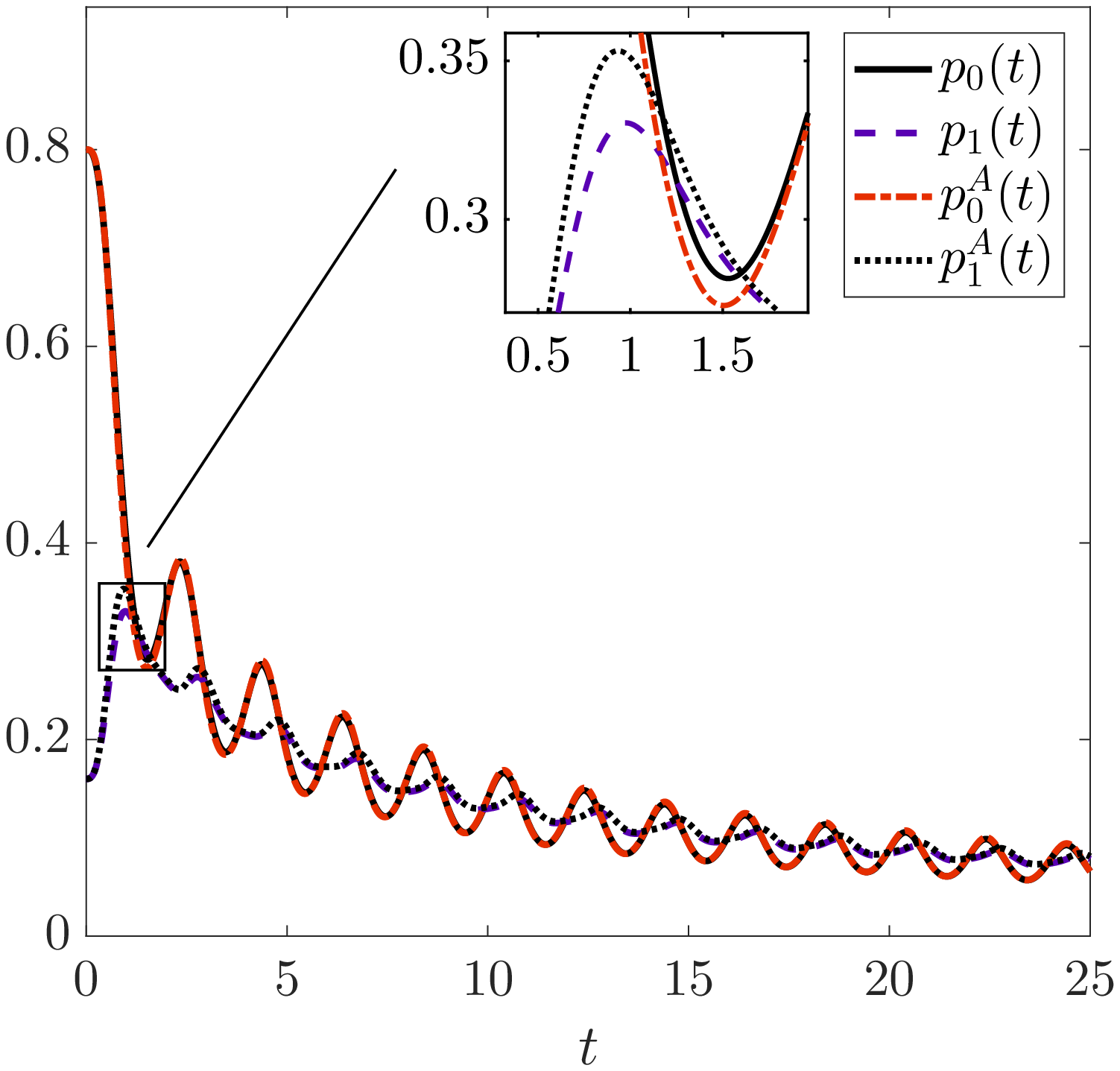}
}
\caption{Comparison of queue-length distribution with exact ODE-system and continuous approximation in the very strong cyclic inflow rate case}
\label{fig:QueueLengthDistVeryStrongCyclic}
\end{figure}
In the very strong cyclic case $\gl_0 = 0.2$ and $\gl_1 = 2$,  we have finite time periods in which the utilization is greater or equal to one and we are in the unstable regime. Nevertheless, visually, the approximations $p_0^A(t)$ and $p_1^A(t)$ are close to the ODE system values, which shows the robustness of the approximation with respect to different utilizations; see figure \ref{fig:QueueLengthDistVeryStrongCyclic}. The errors in table \ref{tab:ErrorTimeDepQueueingModerateCyclic} third column are again of order $10^{-2}$.

\section{Application to Production}
The previous section addressed the formal derivation and numerical validation of the continuous approximation given by equations \eqref{eq:AppQueueModelEq1}-\eqref{eq:AppQueueModelEq3}. In this section, we consider  the $M_t/M_t/1$ queueing model in a production context and we derive measures to evaluate them. The interpretation of this queueing model in the production is as follows: we assume that parts arrive (from, e.g., orders) randomly with a mean rate $\lambda(t)$ and  are put into a waiting queue,  if the production is busy (here, one unit) or into the processor in the case of an idle production unit.  The production time is  random, with a mean rate $\mu(t)$, and the products are fed into the processor from the storage using a FIFO rule. The analysis of the queue length and number of parts in the system has been well-established, and approximations of the expected number of parts are known. One example are the so-called \textit{Pointwise Stationary Fluid Flow Approximations} by, e.g., \cite{Rider1976, Wang1996}.
In addition to the length of the queue, the outflow of the production is the most important measure in a production context. 
%
%

 The outflow in $[t,t+\Delta t]$ is denoted by \[Out(t) = \frac{\text{Number of parts leaving in } [t,t+\Delta t]}{\Delta t}\]
for some $t \geq 0$ and $ \Delta t >0$,  which is a $P$-a.s.\ finite random variable. We can compute the expected outflow in $[t,t+\Delta t]$ as
\begin{align*}
\mathbb{E}[Out(t)] &= \frac{1\cdot P(L(t+\Delta t) = L(t)-1)+o(\Delta t)}{\Delta t}\\
	&=\frac{\sum_{k=1}^\infty P(L(t+\Delta t)=k-1|L(t)=k)P(L(t)=k)}{\Delta t} +o(1)\\
	&=\frac{\sum_{k=1}^\infty \mu(t) \Delta tP(L(t)=k)}{\Delta t} +o(1)\\
	&=\mu(t)(1-P(L(t)=0))+o(1)
\end{align*}
using \eqref{eq:KFTransientMM11} and \eqref{eq:KFTransientMM12}. In a natural manner, we define the expected outflow at time $t$ as the limit $\Delta t \to 0$, i.e., 
\begin{align}
\overline{Out}(t) = \mu(t)(1-p_0(t)). \label{eq:QueueingExpectedOutflowDef}
\end{align}
If we compare the latter with \eqref{eq:queueingExpectation}, we can write the change rate of the expected number of parts in the system as
\begin{align}
\frac{d}{dt} \E[L(t)] = \gl(t)-\overline{Out}(t).  \label{eq:QueueingExpectedQueueLength2}
\end{align}

In numerical experiments, it turned out that using \eqref{eq:QueueingExpectedOutflowDef} in \eqref{eq:QueueingExpectedQueueLength2} leads to avoidable numerical errors in the continuous approximation case.
To calculate the expected number of parts in the system in the case of the continuous approximation, we use the following idea:
\begin{align*}
\frac{d}{dt} \E[L(t)] &= \frac{d}{dt}\sum_{k = 0}^\infty k p_k(t)\\
&\approx \frac{d}{dt}\sum_{k = 0}^\infty k p^A_k(t)\\
&= \sum_{k=0}^\infty k \int_k^{k+1} \rho_t(x,t)dx\\
&= \sum_{k=0}^\infty k (-a(t) \rho(k+1,t)+b(t) \rho_x(k+1,t)+a(t) \rho(k,t)-b(t) \rho_x(k,t))\\
&= a(t) \sum_{k=1}^\infty \rho(k,t) -b(t) \sum_{k=1}^\infty \rho_x(k,t).
\end{align*}
If we consider the numerical approximation, we use the centered difference to approximate $\rho_x(k,t)$ for every $k \in \N$.

In \cite{Rider1976}, a simple approximation of the expected number of parts is derived and is considered in \cite{Wang1996} as well. Let us denote by $\overline{L}^K(t)$ the approximation of the expected number of parts at time $t \geq 0$; then, in \cite{Rider1976}, the approximation satisfies the following initial value problem:
\begin{align}
\frac{d}{dt} \overline{L}^K(t) &= \mu(t) e^{-\mu(t) \tilde{T}}\left(\varrho(t)-\frac{\overline{L}^K(t)}{1+\overline{L}^K(t)}\right),\label{eq:QueueRider}\\
\overline{L}^K(0) &= \overline{L}^K_0 \notag
\end{align}
for some initial value $\overline{L}^K_0 \geq 0$, and $\tilde{T}$ is some parameter used to control the transition. As in the examples in \cite{Rider1976}, we set $\tilde{T} = 0$ and $\overline{L}^K_0 =  \frac{\varrho(0)}{1-\varrho(0)}$, which is the expected queue length in the steady state determined by $\varrho(0)$. 

Equation \eqref{eq:QueueRider} is a first order approximation for the expected queue length. In \cite{MasseyPender2013}, a second and third order approximation are introduced and numerically analyzed. The second order approximation also includes the variance and is called \textit{Gaussian variance approximation} (GVA). The third order model governs additionally the skewness of the queue length distribution and is called \textit{Gaussian skewness approximation} (GSA). For explicit formulas and further details we refer to \cite{MasseyPender2013}. In the following, we denote by $\overline{Out}^{GVA}(t)$ and $\overline{Out}^{GSA}(t)$ the approximate outflow by the GVA and GSA, respectively.

\begin{figure}[H]
\subfigure[Expected outflow approximations]{
\includegraphics[width = 0.95\textwidth]{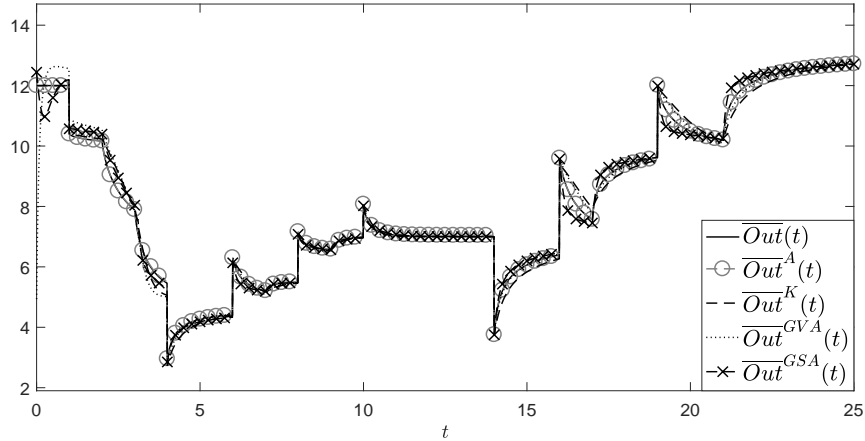}
}
\subfigure[Inflow and processing rates from \cite{Rider1976}]{
\includegraphics[width = 0.4\textwidth]{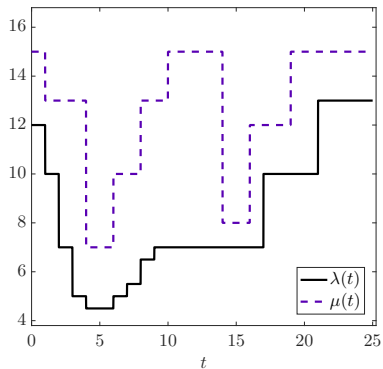}
}
\subfigure[Zoom into expected outflow approximations]{
\includegraphics[width = 0.4\textwidth]{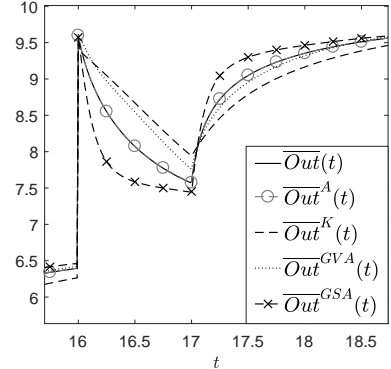}
}
\caption{Comparison of the expected queue length and outflow with inflow and processing rate from \cite{Rider1976}}
\label{fig:ExpectedOutflowKenneth}
\end{figure}

We consider the example in \cite{Rider1976}, where the inflow rate and the processing rate are given by figure \ref{fig:ExpectedOutflowKenneth} (b). For the numerical approximation of our PDE model we used a spatial discretization specified by $\Delta x = 0.01$, $x_0 = 0$ and $x_1 = 200$.

In figure \ref{fig:ExpectedOutflowKenneth} (a), we compare the expected outflow of five approximations:
\begin{enumerate}
\item  $\overline{Out}(t)$ using the ODE system \eqref{eq:KFTransientMM11}-\eqref{eq:KFTransientMM12},  with $N = 1000$, using the closure discussed in \eqref{eq:RestrictODESystem}.
We consider this  the ``exact'' solution. 
\item $\overline{Out}^A(t)$ the continuous approximation \eqref{eq:AppQueueModelEq1}-\eqref{eq:AppQueueModelEq3}.
\item $\overline{Out}^K(t)$ from \eqref{eq:QueueRider}, 
\item $\overline{Out}^{GVA}(t)$ and 
\item $\overline{Out}^{GSA}(t)$, the last two from \cite{MasseyPender2013}.
\end{enumerate}
At a first glance, all expected outflow approximations are quite close. Zooming in figure \ref{fig:ExpectedOutflowKenneth} (c) shows that the continuous approximation $\overline{Out}^A(t)$ coincides with the ODE system result $\overline{Out}(t)$.
The supremum and $L^1$ error between the expected outflow computed with the ODE system and the approximations is shown in table \ref{tab:ErrorKennethInflowAndProcessing}. The high supremum errors for the GVA and GSA are due to the initial instabilities during the first time units.

\begin{table}[htb!]
\centering
\begin{tabular}{c|cccc}
&$\overline{Out}^{A}(t)$&$\overline{Out}^{K}(t)$&$\overline{Out}^{GVA}(t)$&$\overline{Out}^{GSA}(t)$ \\\hline
$\|\cdot\|_\infty$&$0.0457$&$0.6830$&$7.0672$&$1.0314$\\
$\|\cdot\|_{L^1}$&$0.3522$&$3.9450$&$3.8772$&$3.4053$
\end{tabular}
\caption{Error between the approximations and the ODE system result for inflow from \cite{Rider1976}}
\label{tab:ErrorKennethInflowAndProcessing}
\end{table}

\begin{figure}[H]
\subfigure[Expected outflow approximations]{
\includegraphics[width = 0.95\textwidth]{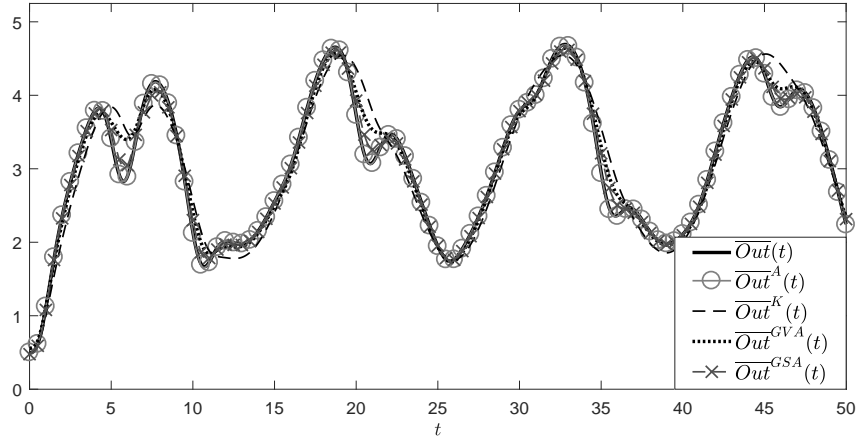}
}
\subfigure[Inflow and processing rates]{
\includegraphics[width = 0.4\textwidth]{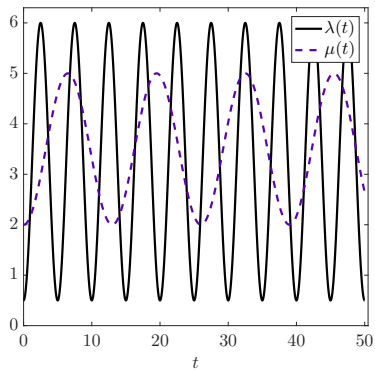}
}
\subfigure[Zoom into expected outflow approximations]{
\includegraphics[width = 0.4\textwidth]{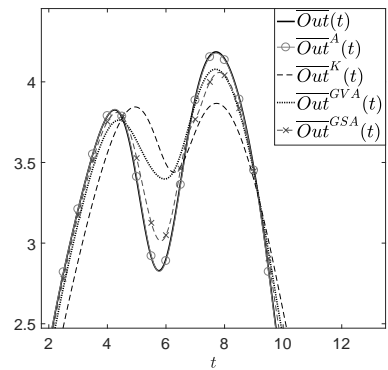}
}
\caption{Comparison of the expected queue length and outflow with cyclic inflow and processing rate}
\label{fig:ExpectedOutflowCycle}
\end{figure}

Finally, we discuss the expected outflow of the queue approximations for a cyclic inflow and processing rate shown in figure \ref{fig:ExpectedOutflowCycle} (b).
There exist over-saturated ($\varrho(t)\geq 1$) and under-saturated ($\varrho(t)<1$) time periods, which imply a strong fluctuation in time. For the expected outflow, we observe in figures \ref{fig:ExpectedOutflowCycle} (a) and (c) that the continuous approximation $\overline{Out}^A(t)$ is again well performing compared to the other approximations. Specifically, the maximal absolute error for the expected outflow is $0.0196$, see table \ref{tab:ErrorCyclicInflow}.
In this example, the simple approximation $\overline{Out}^K(t)$ in \cite{Rider1976} fails to capture the correct dynamics, resulting in an $L^1$ error of $10.5651$ for the expected outflow. Concerning the GVA and GSA, only the GSA  captures the correct dynamic behavior but is not that close to the ODE system as our continuous approximation.

\begin{table}[H]
\centering
\begin{tabular}{c|cccc}
&$\overline{Out}^{A}(t)$&$\overline{Out}^{K}(t)$&$\overline{Out}^{GVA}(t)$&$\overline{Out}^{GSA}(t)$ \\\hline
$\|\cdot\|_\infty$&$0.0196$&$0.9041$&$0.5876$&$0.2901$\\
$\|\cdot\|_{L^1}$&$0.2622$&$10.5651$&$4.8328$&$2.8259$
\end{tabular}
\caption{Error between the approximations and the ODE system result for cyclic inflow}
\label{tab:ErrorCyclicInflow}
\end{table}

\section{Conclusion}
We have derived a continuous approximation of the queue length distributions of an $M_t/M_t/1$ queueing system based on a Fokker-Planck type of partial differential equation under simple assumptions. 

It is instructive to compare our model with the heavy traffic model of queuing theory:
the probability density given by \eqref{eq:AppQueueModelEq1}-\eqref{eq:AppQueueModelEq3} and with the choices \eqref{eq:QueueingKoeffA}-\eqref{eq:QueueingKoeffB} corresponds to a reflected Brownian Motion with drift $\lambda-\mu$ and variance $ \frac{2(\mu(t)-\lambda(t))}{\ln(\mu(t))-\ln(\lambda(t))}$.
The diffusion limit of single-station queues from the literature, see, e.g.,\ \cite{chen2013fundamentals}, leads to a queue-length approximation with drift $\lambda-\mu$ and variance $\lambda+\min\{\lambda,\mu\}$ in the case of exponentially distributed inter-arrival and service times.
In the heavy traffic limit $\varrho \nearrow 1$ both approximations coincide, whereas for $\varrho<1$ our approximation, see \eqref{eq:QueueingKoeffB}, leads to  a higher variance in the system.

We have shown in various numerical examples that our model approximates  the original distribution very well and thus  provides a new approach to  the forward problem of  production planning. In addition  the Fokker-Planck equation  for an initial step at $t=0$ can be solved analytically, allowing the study of  solutions for the transient behavior of the continuous approximation as well as the original queueing system. We introduced an appropriate numerical discretization scheme for the approximate model to study the fully transient cases and compare the solutions of the PDE to the  solution of  a truncated ODE-system for the queue-length distributions.
The PDE approximation shows excellent agreement with a truncation of the queue-length distributions at a 1000 modes. 

Deriving the output from the Fokker-Planck model relates the model production systems. A comparison of the expected outflow of our model, the true expected outflow based on 1000 ODES  and and another well established approximation from the literature shows that our model significantly improves on the literature for all cases considered. 

There are several open avenues for future work:
\begin{itemize}

\item 
In \cite{KoGautam2013}, a diffusion approximation for multi-server queues is introduced using an adjusted fluid and diffusion limit leading to  a system of ODEs. 
We will develop our PDE approximation for multi-server queues and compare it to this and the third order model of   \cite{MasseyPender2013}.

\item 
The improvement of the continuous approximation model over heavy traffic models is essentially due to the fit of the approximation to a known 
stationary distribution of the $M/M/1$ model. There are more complicated queueing networks that have stationary distributions that would be 
obvious candidates for the development of such continuum models. 
In addition  fluid and diffusion limit equations have been introduced for networks \cite{Mandelbaum1998} and can be compared to PDE approximations of such networks. 

\item Having a PDE model that is a good approximation to a queueing system and that even has explicit solutions for some 
relevant cases  allows us to use a wealth of  PDE methods to study this and more complicated  queueing systems. Specifically, 
the production planning problem now becomes an optimal control problem that can be solved via variational methods \cite{la2010control}.

\item Other discrete systems, in particular multi-agent systems, often show a mixture of transport and queuing features that are not resolved well in 
time-dependent and transient cases. We expect that with this approach we can derive better models for traffic and pedestrian flows \cite{Bellomo2011,ChristianiPiccoliTosin2014,GaravelloPiccoli}.
\end{itemize}

\section*{Acknowledgments}
We would like to thank David Kaspar for pointing out the relationship to the heavy traffic limit of queueing theory for us. 
D.A.\ gratefully acknowledges support through NSF grant DMS-1515592 and travel support through the KI-Net grant, NSF RNMS grant No.\ 1107291.
S.G.\ and S.K.\ gratefully thank the BMBF project ENets (05M18VMA) and the DAAD
project ``Stochastic dynamics for complex networks and systems'' (Project-ID 57444394) for the financial support.
\bibliographystyle{siam}
\bibliography{references}
\end{document}